\documentclass[11pt]{amsart}
\usepackage{amscd}
\usepackage{amssymb}
\usepackage{graphicx}
\newcounter{mtheorem}

\newtheorem{theorem}{Theorem}[section]
\newtheorem{lemma}[theorem]{Lemma}
\newtheorem{prop}[theorem]{Proposition}
\newtheorem{corollary}[theorem]{Corollary}

\theoremstyle{definition}
\newtheorem{definition}[theorem]{Definition}
\newtheorem{example}[theorem]{Example}
\newtheorem{decomp}{Decomposition}
\theoremstyle{remark}
\newtheorem{remark}[theorem]{Remark}

\numberwithin{equation}{section}
\newcommand{\E}{\mathcal{E}}
\newcommand{\cone}{\mathcal{C}}

\newcommand{\unitary}[1]{\textrm{U({#1})}}
\newcommand{\sunitary}[1]{\textrm{SU({#1})}}

\newcommand{\C}{\mathbb{C}}
\newcommand{\R}{\mathbb{R}}

\newcommand{\N}{\mathbb{N}}

\newcommand{\Sph}{\mathbb{S}}

\newcommand{\Imag}{\operatorname{Im}}
\newcommand{\Real}{\operatorname{Re}}

\newcommand{\tnabla}{{\widetilde{\nabla}}}
\newcommand{\tg}{{\tilde{g}}}

\newcommand{\tGamma}{{\widetilde{\Gamma}}}
\title[SL conifolds, I]{Special Lagrangian conifolds, I: Moduli spaces}
\author[T.~Pacini]{Tommaso~Pacini}
\address{Scuola Normale Superiore, Pisa} \email{tommaso.pacini@sns.it}
\date{\today}
\begin{document}
\begin{abstract}
We discuss the deformation theory of special Lagrangian (SL) conifolds in $\C^m$. Conifolds are a key ingredient in the compactification problem for moduli spaces of compact SLs in Calabi-Yau manifolds. This category
allows for the simultaneous presence of conical singularities and of non-compact, asymptotically conical, ends. 

Our main theorem is the natural next step in the chain of results
initiated by McLean \cite{mclean} and continued by the author \cite{pacini:defs}
and Joyce \cite{joyce:II}. We emphasize a unifying framework for studying the
various 
cases and discuss analogies and differences between them. This paper also lays down the geometric foundations for
our paper \cite{pacini:slgluing} concerning gluing constructions for SL conifolds in $\C^m$. 
\end{abstract}
\maketitle
\tableofcontents
\section{Introduction}\label{s:intro}
Let $M$ be a Calabi-Yau (CY) manifold. Roughly speaking, a submanifold
$L\subset M$ is \textit{special Lagrangian} (SL) if it is both minimal and
Lagrangian with respect to the ambient Riemannian and symplectic structures. 

From the point of view of Riemannian Geometry it is of course natural to focus
on the minimality condition. It turns out that SLs are automatically
volume-minimizing in their homology class. In fact, this was Harvey and Lawson's
main motivation for defining and studying SLs within the general context of
Calibrated Geometry \cite{harveylawson}. This is still the most common point of
view on SLs and leads to emphasizing the role of analytic and Geometric Measure
Theory techniques. It also provides a connection with various classical problems
in Analysis such as the Plateau problem and the study of area-minimizing cones.
In many ways it is the point of view adopted here.

From the point of view of Symplectic Geometry it is instead natural to focus on
the Lagrangian condition. Specifically, SLs are examples of Maslov-zero
Lagrangian submanifolds. This leads to emphasizing the role of Symplectic
Topology techniques, both classical (such as the h-principle and moment maps)
and contemporary (such as Floer homology). An early instance of this point of
view is the work of Audin \cite{audin}; it also permeates the paper
\cite{haskinspacini} by Haskins and the author.

Given this richness of ingredients it is perhaps not surprising that SLs are
conjectured to play an important role in Mirror Symmetry \cite{kontsevich},
\cite{syz} and to produce interesting new invariants of CY manifolds
\cite{joyce:3spheres}. Likewise, and more intrinsically, they also tend to
exhibit other nice technical features. In particular it is by now well
understood that SLs often generate smooth, finite-dimensional, moduli spaces.
This SL deformation problem has been studied by a number of authors under
various topological and geometric assumptions. One clear path is the chain of
results initiated by McLean \cite{mclean}, who studied deformations of smooth
compact SLs; continued by the author \cite{pacini:defs} and Marshall
\cite{marshall}, who adapted that set-up to study certain smooth non-compact
(\textit{asymptotically conical}, AC) SLs; and further advanced by Joyce, who
presented analogous results for compact \textit{conically singular} (CS) SLs
\cite{joyce:II}. 

The above three classes of SLs are intimately linked, as follows. One of the main open questions in SL geometry is how to compactify McLean's moduli spaces. This problem is currently one of the biggest obstructions to progress on the above conjectures. Roughly speaking, compactifying the moduli space requires adding to it a ``boundary'' containing singular compact SLs. By definition, CS SLs have isolated singularities modelled on SL cones in $\C^m$:  they would be the simplest objects appearing in this boundary. If a CS SL appears in the boundary, it must be a limit of a 1-parameter family of smooth compact SLs. These smooth SLs can be recovered via a gluing construction which desingularizes the CS SL: (i) each singularity of the CS SL defines a SL cone in $\C^m$; (ii) each of these cones must admit a 1-parameter family of SL desingularizations, \textit{i.e.} AC SLs in $\C^m$ converging to the cone as the parameter $t$ tends to $0$; (iii) the family of smooth SLs is obtained by gluing the AC SLs into a neighbourhood of the singularities of the CS SL. This picture is made precise by Joyce's gluing results \cite{joyce:III}, \cite{joyce:IV},
\cite{joyce:V}. Section 8 of \cite{joyce:V} then shows that, in some cases and near the boundary, 
the compactified moduli space can be locally written as a product of moduli
spaces of AC and CS SLs.

The above classes of submanifolds are special cases within the broader
category of \textit{Riemannian conifolds}, which includes manifolds
exhibiting both AC and CS ends. In other words, it allows CS SLs to become non-compact by allowing the presence of AC ends. This is of fundamental importance for the construction of SLs in $\C^m$: it is well-known that $\C^m$ does not admit any compact (smooth or singular) volume-minimizing submanifolds. Cones in $\C^m$ with an isolated singularity at
the origin are the simplest example of conifold: the construction of new examples and the study of their properties is currently one of the most active areas of SL research \cite{harveylawson},
\cite{haskins}, \cite{haskinskapouleas}, \cite{haskinspacini}, \cite{joyce:symmetries}, \cite{ohnita}. Conifolds provide the appropriate framework in which to extend all the above research. In particular, they might also substitute AC SLs in Joyce's gluing results: one could try to cut out a conical singularity of the CS SL and replace it with a different singular conifold, thus jumping from one area of the boundary of the compactified moduli space, containing certain CS SLs, to another. 

The paper at hand is Part I of a multi-step project aiming to set up a general theory of SL conifolds. Two other papers related to this project are currently available: \cite{pacini:weighted}, \cite{pacini:slgluing}. Further work is in progress. The goal of this paper is to provide a general deformation theory of SL conifolds in $\C^m$. The best set-up for the SL deformation
problem is the one provided by Joyce \cite{joyce:II}. It is based on his
Lagrangian neighbourhood and regularity theorems \cite{joyce:I}. Joyce's
framework has two benefits: (i) it simplifies the Analysis via a reduction from
the semi-elliptic operator $d\oplus d^*$ on 1-forms to the
elliptic Laplace operator on functions, (ii) it nicely emphasizes the separate
contributions to the dimension of $\mathcal{M}_L$ coming from the topological
and from the analytic components. After presenting our main result Theorem
\ref{thm:accssl} concerning moduli spaces of CS/AC SL submanifolds in $\C^m$, we thus sketch proofs of the previously-known results emphasizing this point of view. 
In this sense, this paper also serves the purpose of surveying and unifying those results. More importantly, it lays down the geometric foundations for \cite{pacini:slgluing}; the analytic foundations are provided by \cite{pacini:weighted}.

We now summarize the contents of this paper. Section \ref{s:geometry_review}
introduces and studies the category of $m$-dimensional Riemannian conifolds. In particular, Section \ref{s:reviewlaplace} summarizes useful facts concerning the Laplace operator on conifolds while Sections \ref{ss:closedforms} and \ref{ss:closedforms_bis} contain an investigation into the structure of various spaces of closed 1-forms on these
manifolds. This is a fundamental ingredient in the Lagrangian and SL deformation
theory. The corresponding notion of ``subconifolds'' is presented in Section
\ref{s:lagconifolds}, which defines the concept of \textit{Lagrangian conifold}.
Section \ref{s:lagdefs} studies the (infinite-dimensional) deformation theory of Lagrangian conifolds: this relies on 
Joyce's Lagrangian neighbourhood theorems coupled with the material of Section
\ref{ss:closedforms}. After presenting the necessary definitions in Section \ref{s:slgeometry},
the analogous framework for deforming SL conifolds is developed in Section
\ref{s:setup}. The SL deformation theory is
then completed in Section \ref{s:moduli}. Section \ref{ss:comparison} reviews previous results concerning SL moduli spaces, providing a panoramic overview of SL deformation theory. 

To conclude, we should again emphasize that the proof of Theorem \ref{thm:accssl} rests upon three rather delicate and technical ingredients: (i) carefully chosen Lagrangian neighbourhood theorems, (ii) Joyce's SL regularity results and (iii) the theory of weighted Sobolev spaces and elliptic operators on conifolds. In the interest of brevity, in this paper we have kept the presentation of these results to a bare minimum but anyone wishing to do further work in this field will need a deeper understanding of this material. Concerning (i), we thus refer the reader to an expanded version of this paper, available online \cite{pacini:sldefs}. Concerning (ii), we refer the reader to \cite{joyce:I}. Finally, our paper \cite{pacini:weighted} provides full details of the necessary analytic machinery.

\ 

\textbf{Important remark:} To simplify certain arguments, throughout this paper
we assume $m\geq 3$.


\section{Geometry and analysis of conifolds}\label{s:geometry_review}

We introduce here the categories of differentiable and Riemannian manifolds
mainly relevant to this paper, referring to \cite{pacini:weighted} for further
details. Following \cite{joyce:I}, however, we introduce a small variation of
the notion of ``conically singular" manifolds:  presenting them in terms of the
compactification $\bar{L}$ will allow us to keep track of the singular points
$x_i$. This plays no role in this section but in Section \ref{s:lagdefs} it will
become very useful.

\begin{definition}\label{def:manifold_ends}
Let $L^m$ be a smooth manifold. We say $L$ is a \textit{manifold with ends} if
it satisfies the following conditions:
\begin{enumerate}
\item We are given a compact subset $K\subset L$ such that $S:=L\setminus K$ has
a finite number of connected components $S_1,\dots,S_e$, \textit{i.e.}
$S=\amalg_{i=1}^e S_i$.
\item For each $S_i$ we are given a connected ($m-1$)-dimensional compact
manifold $\Sigma_i$ without boundary. 
\item There exist diffeomorphisms $\phi_i:\Sigma_i\times [1,\infty)\rightarrow
\overline{S_i}$.
\end{enumerate}
We then call the components $S_i$ the \textit{ends} of $L$ and the manifolds
$\Sigma_i$ the \textit{links} of $L$. We denote by $S$ the union of the ends and
by $\Sigma$ the union of the links of $L$. 
\end{definition}

\begin{definition}\label{def:metrics_ends}
Let L be a manifold with ends. Let $g$ be a Riemannian metric on $L$. Choose an
end $S_i$ with corresponding link $\Sigma_i$.

We say that $S_i$ is a \textit{conically singular} (CS) end if the following
conditions hold:
\begin{enumerate}
\item $\Sigma_i$ is endowed with a Riemannian metric $g_i'$.

We then let $(\theta,r)$ denote the generic point on the product manifold
$C_i:=\Sigma_i\times (0,\infty)$ and $\tg_i:=dr^2+r^2g_i'$ denote the
corresponding \textit{conical metric} on $C_i$.
\item There exist a constant $\nu_i>0$ and a diffeomorphism
$\phi_i:\Sigma_i\times (0,\epsilon]\rightarrow \overline{S_i}$ such that, as $r\rightarrow
0$ and for all $k\geq 0$,
$$|\tnabla^k(\phi_i^*g-\tg_i)|_{\tg_i}=O(r^{\nu_i-k}),$$
where $\tnabla$ is the Levi-Civita connection on $C_i$ defined by $\tg_i$. 
\end{enumerate}
We say that $S_i$ is an \textit{asymptotically conical} (AC) end if the
following conditions hold:
\begin{enumerate}
\item $\Sigma_i$ is endowed with a Riemannian metric $g_i'$.

We again let $(\theta,r)$ denote the generic point on the product manifold
$C_i:=\Sigma_i\times (0,\infty)$ and $\tg_i:=dr^2+r^2g_i'$ denote the
corresponding conical metric on $C_i$.
\item There exist a constant $\nu_i<0$ and a diffeomorphism
$\phi_i:\Sigma_i\times [R,\infty)\rightarrow \overline{S_i}$ such that, as $r\rightarrow
\infty$ and for all $k\geq 0$,
$$|\tnabla^k(\phi_i^*g-\tg_i)|_{\tg_i}=O(r^{\nu_i-k}),$$
where $\tnabla$ is the Levi-Civita connection on $C_i$ defined by $\tg_i$.
\end{enumerate}
In either of the above situations we call $\nu_i$ the \textit{convergence rate}
of $S_i$.
\end{definition}

We refer to \cite{pacini:weighted} Section 6 for a better understanding of the asymptotic
conditions introduced in Definition \ref{def:metrics_ends}.

\begin{definition}\label{def:cs_manifold}
Let $(\bar{L},d)$ be a metric space. $\bar{L}$ is a \textit{Riemannian manifold
with conical singularities} (CS manifold) if it satisfies the following
conditions.
\begin{enumerate}
\item We are given a finite number of points $\{x_1,\dots,x_e\}\in \bar{L}$ such that $L:=\bar{L}\setminus\{x_1,\dots,x_e\}$ has the
structure of a smooth $m$-dimensional manifold with $e$ ends. 

More specifically, we assume given $\epsilon\in (0,1)$ such that any pair of
distinct points satisfies $d(x_i,x_j)>2\epsilon$. Set $S_i:=\{x\in L:
0<d(x,x_i)<\epsilon\}$. We then assume that $S_i$ are the ends of $L$ with
respect to some given connected links $\Sigma_i$. 
\item We are given a Riemannian metric $g$ on $L$ inducing the distance $d$. 
\item With respect to $g$, each end $S_i$ is CS in the sense of Definition
\ref{def:metrics_ends}.
\end{enumerate}
It follows from our definition that any CS manifold $\bar{L}$ is compact. We
will often not distinguish between $\bar{L}$ and $L$, but notice that $(L,g)$ is
neither compact nor complete. We call $x_i$ the \textit{singularities} of
$\bar{L}$. 
\end{definition}

\begin{definition}\label{def:ac_manifold}
Let $(L,g)$ be a Riemannian manifold. $L$ is a \textit{Riemannian manifold with
asymptotically conical ends} (AC manifold) if it satisfies the following
conditions.
\begin{enumerate}
\item $L$ is a smooth manifold with $e$ ends $S_i$ and connected links
$\Sigma_i$.
\item Each end $S_i$ is AC in the sense of Definition \ref{def:metrics_ends}.
\end{enumerate}
\end{definition}
One can check that AC manifolds are non-compact but complete.

\begin{definition} \label{def:accs_manifold}
Let $(\bar{L},d)$ be a metric space. We say that $\bar{L}$ is a
\textit{Riemannian CS/AC manifold} if it satisfies the following conditions. 
\begin{enumerate}
\item We are given a finite number of points $\{x_1,\ldots,x_s\}$ and a number $l$ such that
$L:=\bar{L}\setminus\{x_1,\dots,x_s\}$ has the structure of a smooth
$m$-dimensional manifold with $s+l$ ends. 
\item We are given a metric $g$ on $L$ inducing the distance $d$.
\item With respect to $g$, neighbourhoods of the points $x_i$ have the structure
of CS ends in the sense of Definition \ref{def:metrics_ends}. These are the
``small" ends. We also assume that the remaining ends are ``large",
\textit{i.e.} they have the structure of AC ends in the sense of Definition
\ref{def:metrics_ends}.
\end{enumerate}
We will denote the union of the CS links (respectively, of the CS ends) by
$\Sigma_0$ (respectively, $S_0$) and those corresponding to the AC links and
ends by $\Sigma_\infty$, $S_\infty$. 
\end{definition}

\begin{definition} \label{def:conifold}
We use the generic term \textit{conifold} to indicate any CS, AC or CS/AC
manifold. If $(L,g)$ is a conifold and $C:=\amalg C_i$ is the union of the
corresponding cones as in Definition \ref{def:metrics_ends}, endowed with the
induced metric $\tg$, we say that $(L,g)$ is \textit{asymptotic} to $(C,\tg)$.
\end{definition}

\begin{remark}\label{rem:nondistinct_sings}
 If we think of $\bar{L}$ as a generic compactification of the manifold with ends $L$, we should allow several CS ends to become connected by the addition of a single singular point. In this section this would however constrast with our assumption that our links are connected, which we adopt to simplify notation. Actually in this section this issue is not of particular interest. In geometric applications it becomes more relevant when dealing with immersed conifolds, as in Section \ref{s:lagconifolds}, but there it is easily solved: by definition an immersion is allowed to identify points, so provided we do not explicitly request that the image points $p_i$ be distinct, it is no problem to assume that the $x_i$ are initially distinct.
\end{remark}

Cones in $\R^n$ are of course the archetype of CS/AC manifold, as follows. 
 
\begin{definition} \label{def:cone}
A subset $\bar{\mathcal{C}}\subseteq\R^n$ is a \textit{cone} if it is invariant
under dilations of $\R^n$, \textit{i.e.} if $t\cdot \bar{\mathcal{C}}\subseteq
\bar{\mathcal{C}}$, for all $t\geq 0$. It is uniquely identified by its
\textit{link} $\Sigma:=\bar{\mathcal{C}}\bigcap \Sph^{n-1}$. We will set
$\mathcal{C}:=\bar{\mathcal{C}}\setminus 0$. The cone is \textit{regular} if
$\Sigma$ is smooth. From now on we will always assume this.

Let $g'$ denote the induced metric on $\Sigma$. Then $\mathcal{C}$ with its
induced metric is isometric to $\Sigma\times (0,\infty)$ with the conical metric
$\tg:=dr^2+r^2g'$. In particular $\bar{\mathcal{C}}$ is a CS/AC manifold; it has
as many AC and CS ends as the number of connected components $\Sigma_i$ of
$\Sigma$. Each $\Sigma_i$ thus defines a singular point $x_i$ but these singular
points are not distinct: they all coincide with the origin. Notice that $\Sigma$
is a subsphere $\Sph^{m-1}\subseteq \Sph^{n-1}$ iff $\bar{\mathcal{C}}$ is an
$m$-plane in $\R^n$.

This example illustrates clearly the issue mentioned in Remark \ref{rem:nondistinct_sings}. Strictly speaking, given our definitions, it would be preferable to think of $\mathcal{C}$ as an immersed copy of the abstract manifold $C:=\amalg_{i=1}^s\Sigma_i\times (0,\infty)$. $\bar{C}$ would be obtained by adding one point to each component of $C$, and the immersion would identify these points by mapping them to $0\in\R^n$.
\end{definition}

Let $E$ be a vector bundle over $(L,g)$. Assume $E$ is endowed with a metric and
metric connection $\nabla$: we say that $(E,\nabla)$ is a \textit{metric pair}.
In later sections $E$ will usually be a bundle of differential forms $\Lambda^r$
on $L$, endowed with the metric and Levi-Civita connection induced from $g$. We
can define two types of Banach spaces of sections of $E$, 
referring to \cite{pacini:weighted} for further details regarding the structure
and properties of these spaces.

Regarding notation, given a vector
$\boldsymbol{\beta}=(\beta_1,\dots,\beta_e)\in \R^e$ and $j\in\N$ we set
$\boldsymbol{\beta}+j:=(\beta_1+j,\dots,\beta_e+j)$. We write
$\boldsymbol{\beta}\geq\boldsymbol{\beta}'$ iff $\beta_i\geq\beta_i'$.

\begin{definition}\label{def:csac_sectionspaces}
Let $(L,g)$ be a conifold with $e$ ends. We say that a smooth function
$\rho:L\rightarrow (0,\infty)$ is a \textit{radius function} if $\rho(x)\equiv
r$ on each end, where up to identifications $r$ is the variable introduced in Definition \ref{def:metrics_ends}. Given any vector
$\boldsymbol{\beta}=(\beta_1,\dots,\beta_{e})\in\R^{e}$, choose a function
$\boldsymbol{\beta}:L\rightarrow\R$ which, on each end $S_i$, restricts to the constant
$\beta_i$. 

Given any metric pair $(E,\nabla)$, the \textit{weighted Sobolev spaces} are
defined by
\begin{equation}\label{eq:weighted_sob}
W^p_{k;\boldsymbol{\beta}}(E):=\mbox{Banach space completion of the space
}\{\sigma\in C^\infty(E):\|\sigma\|_{W^p_{k;\boldsymbol{\beta}}}<\infty\},
\end{equation}
where we use the norm
$\|\sigma\|_{W^p_{k;\boldsymbol{\beta}}}:=(\Sigma_{j=0}^k\int_L|\rho^{
-\boldsymbol{\beta}+j}\nabla^j\sigma|^p\rho^{-m}\,\mbox{vol}_g)^{1/p}$. 

The \textit{weighted spaces of $C^k$ sections} are defined by
\begin{equation}\label{eq:weighted_C^k}
C^k_{\boldsymbol{\beta}}(E):=\{\sigma\in C^k(E):
\|\sigma\|_{C^k_{\boldsymbol{\beta}}}<\infty\},
\end{equation}
where we use the norm $\|\sigma\|_{C^k_{\boldsymbol{\beta}}}:=\sum_{j=0}^k
\mbox{sup}_{x\in L}|\rho^{-\boldsymbol{\beta}+j}\nabla^j\sigma|$. Equivalently,
$C^k_{\boldsymbol{\beta}}(E)$ is the space of sections $\sigma\in C^k(E)$ such
that $|\nabla^j \sigma|=O(r^{\boldsymbol{\beta}-j})$ as $r\rightarrow 0$
(respectively, $r\rightarrow\infty$) along each CS (respectively, AC) end. These
are also Banach spaces.

To conclude, the  \textit{weighted space of smooth sections} is defined by
\begin{equation*}
C^\infty_{\boldsymbol{\beta}}(E):=\bigcap_{k\geq 0} C^k_{\boldsymbol{\beta}}(E).
\end{equation*}
Equivalently, this is the space of smooth sections such that $|\nabla^j
\sigma|=O(\rho^{\boldsymbol{\beta}-j})$ for all $j\geq 0$. This space has a natural
Fr\'echet structure. 

When $E$ is the trivial $\R$ bundle over $L$ we obtain weighted spaces of
functions on $L$. We usually denote these by $W^p_{k,\boldsymbol{\beta}}(L)$ and
$C^k_{\boldsymbol{\beta}}(L)$. In the case of a CS/AC manifold we will often
separate the CS and AC weights, writing
$\boldsymbol{\beta}=(\boldsymbol{\mu},\boldsymbol{\lambda})$ for some
$\boldsymbol{\mu}\in \R^s$ and some $\boldsymbol{\lambda}\in \R^l$. We then
write $C^k_{(\boldsymbol{\mu},\boldsymbol{\lambda})}(E)$ and
$W^p_{k,(\boldsymbol{\mu},\boldsymbol{\lambda})}(E)$.
\end{definition}

For these spaces one can prove the validity of the following weighted version of
the Sobolev Embedding Theorems, cf. \cite{pacini:weighted} Corollary 6.8.

\begin{theorem}\label{thm:embedding}
Let $(L,g)$ be an AC manifold. Let $(E,\nabla)$ be a metric pair over $L$.
Assume $k\geq 0$, $l\in\{1,2,\dots\}$ and $p\geq 1$. Set
$p^*_l:=\frac{mp}{m-lp}$. Then, for all
$\boldsymbol{\beta}'\geq\boldsymbol{\beta}$,
\begin{enumerate}
\item If $lp<m$ then there exists a continuous embedding
$W^p_{k+l,\boldsymbol{\beta}}(E)\hookrightarrow
W^{p^*_l}_{k,\boldsymbol{\beta}'}(E)$.
\item If $lp=m$ then, for all $q\in [p,\infty)$, there exist continuous
embeddings $W^p_{k+l,\boldsymbol{\beta}}(E)\hookrightarrow
W^q_{k,\boldsymbol{\beta}'}(E)$.
\item If $lp>m$ then there exists a continuous embedding
$W^p_{k+l,\boldsymbol{\beta}}(E)\hookrightarrow C^k_{\boldsymbol{\beta}'}(E)$. 
\end{enumerate}
Furthermore, assume $kp>m$. Then the corresponding weighted
Sobolev spaces are closed under multiplication, in the following sense. For
any $\boldsymbol{\beta}_1$ and $\boldsymbol{\beta_2}$ there exists $C>0$ such
that, for all $u\in W^p_{k,\boldsymbol{\beta_1}}$ and $v\in
W^p_{k,\boldsymbol{\beta_2}}$,
\begin{equation*}
\|uv\|_{W^p_{k,\boldsymbol{\beta_1}+\boldsymbol{\beta_2}}}\leq
C\|u\|_{W^p_{k,\boldsymbol{\beta_1}}}\|v\|_{W^p_{k,\boldsymbol{\beta_2}}}.
\end{equation*}
Let $(L,g)$ be a CS manifold. Then the same conclusions hold for all
$\boldsymbol{\beta}'\leq\boldsymbol{\beta}$.

Let $(L,g)$ be a CS/AC manifold. Then, setting
$\boldsymbol{\beta}=(\boldsymbol{\mu},\boldsymbol{\lambda})$, the same
conclusions hold for $\boldsymbol{\mu}'\leq\boldsymbol{\mu}$ on the CS ends and
$\boldsymbol{\lambda}'\geq\boldsymbol{\lambda}$ on the AC ends.
\end{theorem}


\subsection{Review of the Laplace operator on conifolds}\label{s:reviewlaplace}

We now summarize some analytic results concerning the Laplace operator on
conifolds, referring to \cite{pacini:weighted} Section 9 for further details and
references.

\begin{definition}\label{def:exceptionalweights}
Let $(\Sigma,g')$ be a compact Riemannian manifold. Consider the cone
$C:=\Sigma\times (0,\infty)$ endowed with the conical metric
$\tilde{g}:=dr^2+r^2g'$. Let $\Delta_{\tilde{g}}$ denote the corresponding
Laplace operator acting on functions.

For each component $(\Sigma_j,g_j')$ of $(\Sigma,g')$ and each $\gamma\in\R$,
consider the space of homogeneous harmonic functions 
\begin{equation}\label{eq:ac_harmonicter}
V^j_{\gamma}:=\{r^\gamma\sigma(\theta):
\Delta_{\tilde{g}}(r^{\gamma}\sigma)=0\}.
\end{equation}
Set $m^j(\gamma):=\mbox{dim}(V^j_\gamma)$. One can show that $m^j_{\gamma}>0$
iff $\gamma$ satisfies the equation
\begin{equation}\label{eq:exceptionalforlaplacian}
\gamma=\frac{(2-m)\pm\sqrt{(2-m)^2+4e_n^j}}{2},
\end{equation}
for some eigenvalue $e_n^j$ of $\Delta_{g_j'}$ on $\Sigma_j$. Given any weight
$\boldsymbol{\gamma}\in \R^e$, we now set $m(\boldsymbol{\gamma}):=\sum_{j=1}^e
m^j(\gamma_j)$. Let $\mathcal{D}\subseteq\R^e$ denote the set of weights
$\boldsymbol{\gamma}$ for which $m(\boldsymbol{\gamma})>0$. We call these the
\textit{exceptional weights} of $\Delta_{\tg}$.
\end{definition}

Let $(L,g)$ be a conifold. Assume $(L,g)$ is asymptotic to a cone $(C,\tg)$ in
the sense of Definition \ref{def:conifold}. Roughly speaking, the fact that $g$
is asymptotic to $\tilde{g}$ in the sense of Definition \ref{def:metrics_ends}
implies that the Laplace operator $\Delta_g$ is ``asymptotic" to $\Delta_{\tg}$.
Applying Definition \ref{def:exceptionalweights} to $C$ defines weights
$\mathcal{D}\subseteq\R^e$: we call these the \textit{exceptional weights} of
$\Delta_g$. This terminology is due to the following result.

\begin{theorem}\label{thm:laplaceresults}
Let $(L,g)$ be a conifold with $e$ ends. Let $\mathcal{D}$ denote
the exceptional weights of $\Delta_g$. Then $\mathcal{D}$ is a discrete subset of $\R^e$ and the Laplace operator 
\begin{equation*}
\Delta_g:W^p_{k,\boldsymbol{\beta}}(L)\rightarrow
W^p_{k-2,\boldsymbol{\beta}-2}(L)
\end{equation*}
is Fredholm iff $\boldsymbol{\beta}\notin \mathcal{D}$.
\end{theorem}

The above theorem, coupled with a ``change of index formula", leads to the
following conclusion, cf. \cite{pacini:weighted} Section 10.

\begin{corollary}\label{cor:laplaceresults}
Let $(L,g)$ be a compact Riemannian manifold. Consider the map
$\Delta_g:W^p_k(L)\rightarrow W^p_{k-2}(L)$. Then
\begin{equation*}
\mbox{Im}(\Delta_g)=\{u\in W^p_{k-2}(L):\int_L u\,\mbox{vol}_g=0\}, \ \
\mbox{Ker}(\Delta_g)=\R.
\end{equation*}
Let $(L,g)$ be an AC manifold. Consider the map
$\Delta_g:W^p_{k,\boldsymbol{\lambda}}(L)\rightarrow
W^p_{k-2,\boldsymbol{\lambda}-2}(L)$.
If $\boldsymbol{\lambda}>2-m$ is non-exceptional then this map is surjective. If
$\boldsymbol{\lambda}<0$ then this map is injective. Equation \ref{eq:exceptionalforlaplacian} shows that the interval $(2-m,0)$ contains no exceptional weights, so for any 
$\boldsymbol{\lambda}\in (2-m,0)$ it is an isomorphism. 

Let $(L,g)$ be a CS manifold with $e$ ends. Consider the map
$\Delta_g:W^p_{k,\boldsymbol{\mu}}(L)\rightarrow
W^p_{k-2,\boldsymbol{\mu}-2}(L)$. If $\boldsymbol{\mu}\in (2-m,0)$ then
\begin{equation*}
\mbox{Im}(\Delta_g)=\{u\in W^p_{k-2,\boldsymbol{\mu}-2}(L): \int_L u\,vol_g=0\},
\ \ Ker(\Delta_g)=\R.
\end{equation*} 
If $\boldsymbol{\mu}>0$ is non-exceptional then this map is injective and 
\begin{equation*}
\mbox{dim(Coker($\Delta_g$))}=e+\sum_{0<\boldsymbol{\gamma}<\boldsymbol{\mu}}
m(\boldsymbol{\gamma}),
\end{equation*}
where $m(\boldsymbol{\gamma})$ is as in Definition \ref{def:exceptionalweights}.

Let $(L,g)$ be a CS/AC manifold with $s$ CS ends and $l$ AC ends. Consider the
map
\begin{equation*}
\Delta_g:W^p_{k,(\boldsymbol{\mu},\boldsymbol{\lambda})}(L)\rightarrow
W^p_{k-2,(\boldsymbol{\mu}-2,\boldsymbol{\lambda}-2)}(L).
\end{equation*}
If $(\boldsymbol{\mu},\boldsymbol{\lambda})\in (2-m,0)$ then this map is an
isomorphism. If $\boldsymbol{\mu}>0$ and $\boldsymbol{\lambda}<0$ are
non-exceptional then this map is injective and 
\begin{equation*}
\mbox{dim(Coker($\Delta_g$))}=s+\sum_{0<\boldsymbol{\gamma}<\boldsymbol{\mu}}
m(\boldsymbol{\gamma}),
\end{equation*}
where $m(\boldsymbol{\gamma})$ is as in Definition \ref{def:exceptionalweights}.
Notice in particular that this dimension depends only on the harmonic functions
on the CS cones.
\end{corollary}


\subsection{Cohomology of manifolds with ends}\label{ss:closedforms}

Let $L$ be a smooth compact manifold or a smooth manifold with ends. 
Then $L$ has topology of
finite type so the first cohomology group
$$H^1(L;\R):=\frac{\{\mbox{Smooth closed $1$-forms on $L$}\}}{d(C^\infty(L))}$$ 
has finite dimension $b^1(L)$. This proves the following statement concerning the
structure of the space of smooth closed $1$-forms.

\begin{decomp}[for compact manifolds or manifolds with
ends]\label{decomp:closedforms} 
Let $L$ be a smooth compact manifold or a smooth manifold with ends. Choose a
finite-dimensional vector space $H$ of closed $1$-forms on $L$ such that the map
\begin{equation}
H\rightarrow H^1(L;\R), \ \ \alpha\mapsto [\alpha]
\end{equation}
is an isomorphism. Then
\begin{equation}\label{eq:closedforms}
\{\mbox{Smooth closed $1$-forms on $L$}\}= H\oplus d(C^\infty(L)). 
\end{equation}
\end{decomp}

We now want to show that in the case of a manifold with ends there exist natural
conditions on the space of 1-forms $H$. 

\begin{definition}\label{def:translationinvariant}
Given a manifold $\Sigma$, set $C:=\Sigma\times (0,\infty)$. Consider the
projection $\pi:\Sigma\times (0,\infty)\rightarrow \Sigma$. A $p$-form $\eta$ on
$C$ is \textit{translation-invariant} if it is of the form $\eta=\pi^*\eta'$,
for some $p$-form $\eta'$ on $\Sigma$.
\end{definition}

\begin{lemma}\label{lemma:formclosed}
Let $L$ be a smooth manifold with ends $S_i$. Let $\alpha$ be a smooth closed
1-form on $L$. Then there exist a smooth closed 1-form $\alpha'$ and a smooth
function $A$ on $L$ such that $\alpha'_{|S_i}$ is translation-invariant and
$\alpha=\alpha'+dA$. 
If furthermore $\alpha$ has compact support then we can choose $\alpha'$ to have
compact support.
\end{lemma}
\begin{proof}
The proof follows the scheme of the Poincar\'e Lemma for de Rham cohomology, cf.
\textit{e.g.} \cite{botttu}. Given any $p$-form $\eta$ on $S_i=\Sigma_i\times
(1,\infty)$, we can write 
$$\eta=\eta_1(\theta,r)+\eta_2(\theta,r)\wedge dr$$
for some $r$-dependent $p$-form $\eta_1$ and ($p-1$)-form $\eta_2$ on $\Sigma$.
Specifically, $\eta_1$ is the restriction of $\eta$ to the cross-sections
$\Sigma_i\times\{r\}$ and $\eta_2:=i_{\partial r} \eta$. For a fixed $R_0>1$ we
then define $(K\eta)(\theta,r):=\int_{R_0}^r\eta_2(\theta,\rho)\,d\rho$. 

Let us apply this to the 1-form obtained by restricting $\alpha$ to $S_i$,
writing
$$\alpha_{|S_i}=\alpha_1(\theta,r)+\alpha_2(\theta,r)\, dr$$
for some $r$-dependent 1-form $\alpha_1$ and function $\alpha_2$ on $\Sigma_i$.
It is then easy to check that
\begin{eqnarray*}
d\alpha_{|S_i} &=& d_\Sigma\alpha_1-(\frac{\partial}{\partial r} \alpha_1)\wedge
dr+(d_\Sigma\alpha_2)\wedge dr,\\
K\alpha_{|S_i} &=&\int_{R_0}^r\alpha_2(\theta,\rho)\,d\rho,\\
d(K\alpha_{|S_i})
&=&\int_{R_0}^rd_\Sigma\alpha_2(\theta,\rho)\,d\rho+\alpha_2(\theta,r)\,dr.
\end{eqnarray*}
From $d\alpha=0$ it follows that $\alpha_1(\theta,R_0)+d(K\alpha)=\alpha_{|S_i}$
and that $\alpha_1(\theta,R_0)$ is closed. Setting
$\alpha'_i:=\alpha_1(\theta,R_0)$ and $A_i:=K\alpha$ we can rewrite this as
$\alpha_{|S_i}=\alpha'_i+dA_i$. Interpolating between the $A_i$ yields a global
smooth function $A$ on $L$ such that $\alpha_{|S_i}=\alpha'_i+dA_{|S_i}$. We can
now define $\alpha':=\alpha-dA$ to obtain the global relationship
$$\alpha=\alpha'+dA.$$
It is clear from this construction that if $\alpha$ has compact support then
(choosing $R_0$ large enough) $\alpha'$ also has compact support.
\end{proof}

Recall that compactly-supported forms give rise to the following theory. Let $L$
be a smooth manifold with ends. We denote by $\Lambda^p_c(L;\R)$ the space of
smooth compactly-supported $p$-forms on $L$ and by $H^p_c(L;\R)$ the
corresponding cohomology groups. Let $\Sigma$ denote the union of the links of
$L$. Notice that $L$ is deformation-equivalent to a compact manifold with
boundary $\Sigma$. Standard algebraic topology (see also \cite{joyce:I} Section
2.4) proves that the inclusion $\Sigma\subset L$ gives rise to a long exact
sequence in cohomology
\begin{equation}\label{eq:cohomsequence}
0\rightarrow H^0(L;\R)\rightarrow H^0(\Sigma;\R)\stackrel{\delta}{\rightarrow}
H^1_c(L;\R)\stackrel{\gamma}{\rightarrow} H^1(L;\R)\stackrel{\rho}{\rightarrow}
H^1(\Sigma;\R)\rightarrow\dots.
\end{equation}
Here, $\gamma$ is induced by the injection $\Lambda^1_c(L;\R) \rightarrow
\Lambda^1(L;\R)$ and $\rho$ is induced by the restriction
$\Lambda^1(L;\R)\rightarrow \Lambda^1(\Sigma;\R)$. We set
$\widetilde{H}^1_c:=\mbox{Im}(\gamma)=\mbox{Ker}(\rho)$. Exactness implies that 
\begin{align}\label{eq:dimexactseq}
\mbox{dim}(\widetilde{H}^1_c) &=
\mbox{dim}(H^1_c(L;\R))-\mbox{dim}(H^0(\Sigma;\R))+\mbox{dim}(H^0(L;\R))\\
\nonumber
&= b^1_c(L)-e+1.
\end{align}
\begin{remark} \label{rem:compactsupport}
The sequence \ref{eq:cohomsequence} shows that 
\begin{equation}\label{eq:exactseq}
H^1_c(L,\R)\simeq \widetilde{H}^1_c\oplus \mbox{Ker}(\gamma) = \widetilde{H}^1_c
\oplus \mbox{Im}(\delta).
\end{equation}
This decomposition can be expressed in words as follows. By definition,
$H^1_c(L;\R)$ is determined by the classes of compactly-supported 1-forms which
are not the differential of a compactly-supported function. Given any such form,
there are two cases: (i) it is not the differential of \textit{any} function, in
which case $\gamma$ maps its class to a non-zero element of $\widetilde{H}^1_c$,
(ii) it is the differential of some function, in which case $\gamma$ maps its
class to zero. However, this function is necessarily constant on the ends of
$L$: these constants can be parametrized via $H^0(\Sigma;\R)$. Notice that the
function is only well-defined up to a constant; likewise, $\mbox{Im}(\delta)$
coincides with $H^0(\Sigma;\R)$ only up to $H^0(L;\R)\simeq\R$.
\end{remark}

Concerning Decomposition \ref{decomp:closedforms}, we can now choose $H$ as follows. For $i=1,\dots,k=\mbox{dim}(\tilde{H}^1_c)$
let $[\alpha_i]$ be a basis of $\widetilde{H}^1_c$. According to Lemma
\ref{lemma:formclosed} we can choose $\alpha_i'$ with compact support such that
[$\alpha_i']=[\alpha_i]$. For $i=1,\dots,N=\mbox{dim}(H^1)$ let $[\alpha_i]$
denote an extension to a basis of $H^1(L;\R)$. Again using Lemma
\ref{lemma:formclosed} we can choose an extension $\alpha_i'$ of
translation-invariant 1-forms such that [$\alpha_i']=[\alpha_i]$. Set
\begin{equation}\label{eq:naturalH}
\widetilde{H}:=\mbox{span}\{\alpha_1',\dots,\alpha_k'\},\ \
H:=\mbox{span}\{\alpha_1',\dots,\alpha_N'\}.
\end{equation}
Then $H$ satisfies the assumptions of Decomposition \ref{decomp:closedforms}.
One advantage of this choice of $H$ is that it reflects the relationship of
$\widetilde{H}^1_c$ to $H^1$. Specifically, if we apply Decomposition
\ref{decomp:closedforms} to $\alpha$ writing $\alpha=\alpha'+dA$ with
$\alpha'\in H$, then $[\alpha]\in \widetilde{H}^1_c$ iff $\alpha'\in
\widetilde{H}$, \textit{i.e.} iff $\alpha'$ has compact support.


\subsection{Cohomology of conifolds}\label{ss:closedforms_bis}

We now want to achieve analogous decompositions for CS and AC manifolds, in
terms of weighted spaces of closed and exact $1$-forms. 

\begin{lemma}\label{l:translationinvariant}
Let $(\Sigma,g')$ be a Riemannian manifold. Let the corresponding cone $C$ have
the conical metric $\tg:=dr^2+r^2g'$. Then any translation-invariant $p$-form
$\eta=\pi^*\eta'$ belongs to the weighted space $C^\infty_{(-p,-p)}(\Lambda^p)$.
For any $\beta>0$, $\eta$ belongs to the smaller weighted space
$C^\infty_{(-p+\beta,-p-\beta)}(\Lambda^p)$ iff $\eta'=0$.
\end{lemma}
\begin{proof} As seen in the proof of Lemma \ref{lemma:formclosed}, the general
$p$-form $\eta$ on $C$ can be written
$\eta=\eta_1(\theta,r)+\eta_2(\theta,r)\wedge dr$. The form is
translation-invariant iff $\eta_1$ is $r$-independent and $\eta_2=0$. In this
case $|\eta|_{\tg}=r^{-p}|\eta_1|_{g'}$ so $|\eta|_{\tg}=O(r^{-p})$ both for
$r\rightarrow 0$ and for $r\rightarrow \infty$. This proves that $\eta\in
C^0_{(-p,-p)}(\Lambda^p)$. To show that $\eta\in C^\infty_{(-p,-p)}(\Lambda^p)$
it is necessary to estimate $|\tnabla^k\eta|_{\tg}$, where $\tnabla$ is the
Levi-Civita connection. This can be done fairly explicitly in terms of
Christoffel symbols. In particular one can choose local coordinates on
$U\subset\Sigma$ defining a local frame $\partial_1,\cdots,\partial_{m-1}$. Set
$\partial_0:=\partial r$, the standard frame on $(0,\infty)$. The Christoffel
symbols for the corresponding frame on $(0,\infty)\times U$ and the metric $\tg$
can then be computed explicitly: for $i,j,k\geq 1$ one finds that
$\tGamma_{i,j}^k$ is bounded, $\tGamma_{i,j}^0=O(r)$,
$\tGamma_{i,0}^k=O(r^{-1})$,
$\tGamma_{0,0}^k=\tGamma_{i,0}^0=\tGamma_{0,0}^0=0$. The Christoffel symbols
defined by $\tg$ for the other tensor bundles depend linearly on these, so they
have the same bounds. Using these calculations one finds that
$|\tnabla^k\eta|_{\tg}=O(r^{-p-k})$, as desired.

It is clear from the proof that $\eta$ satisfies stronger bounds iff it
vanishes.
\end{proof}

\begin{decomp}[for CS or AC manifolds and forms with allowable
growth]\label{decomp:closedforms_growth}
Let $L$ be a CS manifold. Choose a finite-dimensional vector space $H$ of smooth
closed 1-forms on $L$ as in Equation \ref{eq:naturalH}. Then, for any
$\boldsymbol{\beta}<0$,
\begin{equation}\label{eq:closedforms_growth_cs}
\{\mbox{Closed 1-forms on $L$ in $C^\infty_{\boldsymbol{\beta}-1}(\Lambda^1)$}\}
= H\oplus d(C^\infty_{\boldsymbol{\beta}}(L)).
\end{equation}
Analogously, let $L$ be an AC manifold. Choose $H$ as above. Then, for any
$\boldsymbol{\beta}>0$,
\begin{equation}\label{eq:closedforms_growth_ac}
\{\mbox{Closed 1-forms on $L$ in $C^\infty_{\boldsymbol{\beta}-1}(\Lambda^1)$}\}
= H\oplus d(C^\infty_{\boldsymbol{\beta}}(L)).
\end{equation}
\end{decomp}

\begin{proof} 
Consider the CS case. Since $\boldsymbol{\beta}<0$, Lemma
\ref{l:translationinvariant} proves that $H\oplus
d(C^\infty_{\boldsymbol{\beta}}(L))\subseteq \{\mbox{Closed 1-forms in
$C^\infty_{\boldsymbol{\beta}-1}(\Lambda^1)$}\}$. Now choose a closed $\alpha\in
C^\infty_{\boldsymbol{\beta}-1}(\Lambda^1)$. By Decomposition
\ref{decomp:closedforms} we can write $\alpha=\alpha'+dA$, for some $\alpha'\in
H$ and $A\in C^\infty(L)$. Notice that $dA=\alpha-\alpha'\in
C^\infty_{\boldsymbol{\beta}-1}(\Lambda^1)$. By integration, again using the
fact $\boldsymbol{\beta}<0$, we conclude that $A\in
C^\infty_{\boldsymbol{\beta}}(L)$. This proves the opposite inclusion, thus the
identity. The AC case is analogous.
\end{proof}
 
\begin{lemma}\label{l:formexact}
Assume $L$ is a CS manifold. If $\alpha$ is a smooth closed 1-form on $L$
belonging to the space $C^\infty_{\boldsymbol{\beta}-1}(\Lambda^1)$ for some
$\boldsymbol{\beta}>0$ then there exists a smooth closed 1-form $\alpha'$ with
compact support on $L$ and a smooth function $A\in
C^\infty_{\boldsymbol{\beta}}(L)$ such that $\alpha=\alpha'+dA$. 

Assume $L$ is an AC manifold. If $\alpha$ is a smooth closed 1-form on $L$
belonging to the space $C^\infty_{\boldsymbol{\beta}-1}(\Lambda^1)$ for some
$\boldsymbol{\beta}<0$ then there exists a smooth closed 1-form $\alpha'$ with
compact support on $L$ and a smooth function $A\in
C^\infty_{\boldsymbol{\beta}}(L)$ such that $\alpha=\alpha'+dA$. 
\end{lemma}
\begin{proof} The proof is a variation of the proof of Lemma
\ref{lemma:formclosed}, as follows. Consider the AC case. Write
$\alpha_{|S_i}=\alpha_1+\alpha_2\wedge dr$. Define
$K\alpha:=-\int_r^\infty\alpha_2(\theta,\rho)\,d\rho$: this converges because
$\boldsymbol{\beta}<0$. It is simple to check that $d(K\alpha)=\alpha$; in
particular, this shows that $\alpha$ is exact on each end $S_i$. Setting
$A:=K\alpha$ and extending as in Lemma \ref{lemma:formclosed} leads to a global
decomposition $\alpha=\alpha'+dA$ on $L$. By construction $\alpha'$ has compact
support and $A\in C^\infty_{\boldsymbol{\beta}}$. The CS case is analogous, with
$K\alpha:=\int_0^r\alpha_2(\theta,\rho)\,d\rho$.
\end{proof}

\begin{decomp}[for CS or AC manifolds and forms with allowable
decay]\label{decomp:closedforms_decay}
Let $L$ be a CS manifold. Assume $\boldsymbol{\beta}>0$. Choose a
finite-dimensional vector space $H$ of closed 1-forms on $L$ as in Equation
\ref{eq:naturalH}, using $\widetilde{H}_0$ to denote the space $\widetilde{H}$.
For any $i=1,\dots,e$ choose a smooth function $f_i$ on $L$ such that $f_i\equiv
1$ on the end $S_i$ and $f_i\equiv 0$ on the other ends. We can do this in such
a way that $\sum f_i\equiv 1$. Let $E_0$ denote the $e$-dimensional vector space
generated by these functions. By construction $E_0$ contains the constant
functions so $d(E_0)$ has dimension $e-1$. It is simple to check that
$d(E_0)\cap d(C^\infty_{\boldsymbol{\beta}}(L))=\{0\}$. Then
\begin{equation}\label{eq:closedforms_decay_cs}
\{\mbox{Closed 1-forms on $L$ in
$C^\infty_{\boldsymbol{\beta}-1}(\Lambda^1)$}\}=\widetilde{H}_0\oplus
d(E_0)\oplus d(C^\infty_{\boldsymbol{\beta}}(L)).
\end{equation}
Analogously, let $L$ be an AC manifold. Assume $\boldsymbol{\beta}<0$. Choose
spaces as above, this time using the notation $\widetilde{H}_\infty$ and
$E_\infty$. Then
\begin{equation}\label{eq:closedforms_decay_ac}
\{\mbox{Closed 1-forms on $L$ in $C^\infty_{\boldsymbol{\beta}-1}(\Lambda^1)$}\}
= \widetilde{H}_\infty\oplus d(E_\infty)\oplus
d(C^\infty_{\boldsymbol{\beta}}(L)).
\end{equation}
\end{decomp}
\begin{proof} Consider the CS case. The inclusion $\supseteq$ is clear.
Conversely, let $\alpha\in C^\infty_{\boldsymbol{\beta}-1}(\Lambda^1)$ be
closed. Decomposition \ref{decomp:closedforms} allows us to write
$\alpha=\alpha'+dA$, for some uniquely defined $\alpha'\in H$ and some $A\in
C^\infty(L)$, well-defined up to a constant. Lemma \ref{l:formexact} implies
that the cohomology class of $\alpha$ belongs to the space $\widetilde{H}^1_c$,
\textit{i.e.} that $\alpha'\in \widetilde{H}_0$ so it has compact support. This
shows that $dA\in C^\infty_{\boldsymbol{\beta}-1}(\Lambda^1)$. Writing
$A_i:=A_{|S_i}$ we find $dA_i=d_{\Sigma_i}A_i+\frac{\partial A_i}{\partial
r}\,dr$, thus $\frac{\partial A_i}{\partial r}\in
C^\infty_{\boldsymbol{\beta}-1}(L)$. This shows that $\int_0^r\frac{\partial
A_i}{\partial r}\,d\rho\in C^\infty_{\boldsymbol{\beta}}(L)$. This determines
$A_i$ up to a constant $c_i$ on each end. Together with Equation
\ref{eq:exactseq} this proves the claim. The AC case is analogous.
\end{proof}

We now turn to the case of CS/AC manifolds, concentrating on the situations of
most interest to us.

\begin{decomp}[for CS/AC manifolds]\label{decomp:closedforms_accs}
Let $L$ be a CS/AC manifold with $s$ CS ends and $l$ AC ends. As usual we denote
the union of the CS links by $\Sigma_0$ and the union of the AC links by
$\Sigma_\infty$. Choose a finite-dimensional vector space $H$ of closed 1-forms
on $L$ as in Equation \ref{eq:naturalH}, using $\widetilde{H}_{0,\infty}$ to
denote the space $\widetilde{H}$. For any $i=1,\dots,s+l$ choose a function
$f_i$ such that $f_i\equiv 1$ on the end $S_i$ and $f_i\equiv 0$ on the other
ends. We can assume that $\sum f_i\equiv 1$. Let $E_{0,\infty}$ denote the
$(s+l)$-dimensional vector space generated by these functions. Then, for any
$\boldsymbol{\mu}>0$ and $\boldsymbol{\lambda}<0$,
\begin{equation}\label{eq:closedforms_decay_accs}
\{\mbox{Closed 1-forms on $L$ in
$C^\infty_{(\boldsymbol{\mu}-1,\boldsymbol{\lambda}-1)}(\Lambda^1)$}\}=
\widetilde{H}_{0,\infty}\oplus d(E_{0,\infty})\oplus
d(C^\infty_{(\boldsymbol{\mu},\boldsymbol{\lambda})}(L)). 
\end{equation}
Now let $\Lambda^p_{c,\bullet}(L;\R)$ denote the space of $p$-forms on $L$ which
vanish in a neighbourhood of the singularities, with no condition on the large
ends. Let $H^p_{c,\bullet}(L;\R)$ denote the corresponding cohomology groups.
Let $\widetilde{H}^1_{c,\bullet}$ denote the image of the map
$\gamma:H^1_{c,\bullet}(L;\R)\rightarrow H^1(L;\R)$. Choose a finite-dimensional
vector space $\widetilde{H}_{0,\bullet}$ of translation-invariant closed 1-forms
on $L$ with compact support in a neighbourhood of the singularities and such
that the map
\begin{equation}
\widetilde{H}_{0,\bullet}\rightarrow \widetilde{H}^1_{c,\bullet}, \ \
\alpha\mapsto [\alpha]
\end{equation}
is an isomorphism. For any $i=1,\dots,s$ choose a function $f_i$ such that
$f_i\equiv 1$ on the CS end corresponding to the singularity $x_i$ and
$f_i\equiv 0$ on the other ends. Let $E_0$ denote the s-dimensional vector space
generated by these functions. Then, for any $\boldsymbol{\mu}>0$ and
$\boldsymbol{\lambda}>0$,
\begin{equation}\label{eq:closedforms_growth_accs}
\{\mbox{Closed 1-forms on $L$ in
$C^\infty_{(\boldsymbol{\mu}-1,\boldsymbol{\lambda}-1)}(\Lambda^1)$}\}=
\widetilde{H}_{0,\bullet}\oplus d\Big(E_0\oplus
C^\infty_{(\boldsymbol{\mu},\boldsymbol{\lambda})}(L)\Big). 
\end{equation}
\end{decomp}
\proof{} The proof is similar to the proofs of the previous decompositions. It
may however be good to emphasize that, in the case $\boldsymbol{\mu}>0$ and
$\boldsymbol{\lambda}>0$, $d(E_0)\cap
d(C^\infty_{(\boldsymbol{\mu},\boldsymbol{\lambda})}(L))\neq\{0\}$ (it is
one-dimensional). This explains the slightly different statement of
Decomposition \ref{eq:closedforms_growth_accs}.
\endproof

\begin{remark}
The weight $\boldsymbol{\beta}=0$ corresponds to an exceptional case in Lemma
\ref{l:formexact}: integration will generally generate log terms, so we cannot
conclude that $A\in C^\infty_{\boldsymbol{\beta}}$ there. One can analogously
argue that
$C^\infty_{-\boldsymbol{1}}(\Lambda^1)/d(C^\infty_{\boldsymbol{0}}(L))$ is not
finite-dimensional. 

Similar decompositions hold for $k$-forms: in this setting the exceptional case
corresponds to $\boldsymbol{\beta}=k-1$.
\end{remark} 

\begin{remark}
Notice that the above decompositions do not cover all possibilities: for
example, given a CS manifold we could decide to study the space of closed
1-forms in $C^\infty_{\boldsymbol{\beta}-1}(\Lambda^1)$ corresponding to a
weight $\boldsymbol{\beta}=(\beta_1,\dots,\beta_e)$ with some $\beta_i$ positive
and others negative. However, it should be clear from the above discussion how
to use the same ideas to cover any other case of interest. We have restricted
our attention to the cases most relevant to this paper.
\end{remark}

For future reference it is useful to emphasize the topological interpretation of
some of the previous results. The reasons underlying our interest for each case
will become apparent in Section \ref{s:moduli}. 

\begin{corollary}\label{cor:topsummary}
Let $L$ be a smooth compact manifold. Then 
\begin{equation*}
\{\mbox{Closed $1$-forms on $L$}\}\simeq H^1(L;\R)\oplus d(C^\infty(L)). 
\end{equation*}
Let $(L,g)$ be an AC manifold. Then for $\boldsymbol{\beta}<0$
\begin{equation*}
\{\mbox{Closed $1$-forms on $L$ in
$C^\infty_{\boldsymbol{\beta}-1}(\Lambda^1)$}\} \simeq H^1_c(L;\R) \oplus
d(C^\infty_{\boldsymbol{\beta}}(L)), 
\end{equation*}
while for $\boldsymbol{\beta}>0$
\begin{equation*}
\{\mbox{Closed $1$-forms on $L$ in
$C^\infty_{\boldsymbol{\beta}-1}(\Lambda^1)$}\} \simeq H^1(L;\R) \oplus
d(C^\infty_{\boldsymbol{\beta}}(L)).
\end{equation*}
Let $(L,g)$ be a CS manifold with link $\Sigma_0$. Then for
$\boldsymbol{\beta}>0$
\begin{align*}
&\{\mbox{Closed $1$-forms on $L$ in
$C^\infty_{\boldsymbol{\beta}-1}(\Lambda^1)$}\}\\
&\quad\simeq \mbox{Ker}\left(H^1(L)\stackrel{\rho}{\rightarrow}
H^1(\Sigma_0)\right) \oplus d(E_0)\oplus d(C^\infty_{\boldsymbol{\beta}}(L)).
\end{align*}
Let $(L,g)$ be a CS/AC manifold with link $\Sigma=\Sigma_0\amalg\Sigma_\infty$.
Then for $\boldsymbol{\mu}>0$ and $\boldsymbol{\lambda}<0$
\begin{align*}
&\{\mbox{Closed $1$-forms on $L$ in
$C^\infty_{(\boldsymbol{\mu}-1,\boldsymbol{\lambda}-1)}(\Lambda^1)$}\}\\
&\quad\simeq \mbox{Ker}\left(H^1_{\bullet,c}(L)\stackrel{\rho}{\rightarrow}
H^1(\Sigma_0)\right)\oplus d(E_0) \oplus
d(C^\infty_{(\boldsymbol{\mu},\boldsymbol{\lambda})}(L)),
\end{align*}
while for $\boldsymbol{\mu}>0$ and $\boldsymbol{\lambda}>0$
\begin{align*}
&\{\mbox{Closed $1$-forms on $L$ in
$C^\infty_{(\boldsymbol{\mu}-1,\boldsymbol{\lambda}-1)}(\Lambda^1)$}\}\\
&\quad\simeq \mbox{Ker}\left(H^1(L)\stackrel{\rho}{\rightarrow}
H^1(\Sigma_0)\right)\oplus d\left(E_0\oplus
C^\infty_{(\boldsymbol{\mu},\boldsymbol{\lambda})}(L)\right).
\end{align*}
\end{corollary}
\begin{proof}
The compact case coincides with Equation \ref{eq:closedforms}. The AC case with
$\boldsymbol{\beta}<0$ follows from Equation \ref{eq:closedforms_decay_ac} and
Remark \ref{rem:compactsupport}. The AC case with $\boldsymbol{\beta}>0$
coincides with Equation \ref{eq:closedforms_growth_ac}. The CS case coincides
with Equation \ref{eq:closedforms_decay_cs}. 

Let us now focus on the CS/AC case with $\boldsymbol{\lambda}<0$. Using the
notation of Decomposition \ref{decomp:closedforms_accs}, let $E'$ denote a
complement of $E_0\oplus\R$ in $E_{0,\infty}$, \textit{i.e.}
$E_{0,\infty}=E_0\oplus\R\oplus E'$. Notice that the long exact sequence
\ref{eq:cohomsequence} with $\Sigma=\Sigma_0\amalg\Sigma_\infty$ leads to an
identification $H^1_c(L;\R)\simeq \widetilde{H}^1_c(L)\oplus d(E_{0,\infty})$.
One can also set up the ``relative" analogue of Sequence \ref{eq:cohomsequence}
using the inclusion of pairs $(\Sigma_0,\emptyset)\subset (L,\Sigma_\infty)$.
Using notation analogous to that of Decomposition \ref{decomp:closedforms_accs}
this leads to the long exact sequence
\begin{equation*}
0\rightarrow H^0_c(L;\R)\rightarrow H^0_{\bullet,c}(L;\R)\rightarrow
H^0(\Sigma_0;\R)\rightarrow H^1_c(L;\R)\stackrel{\gamma}{\rightarrow}
H^1_{\bullet,c}(L;\R)\stackrel{\rho}{\rightarrow}
H^1(\Sigma_0;\R)\rightarrow\dots
\end{equation*}
Since $H^0_c(L;\R)=0$ and $H^0_{\bullet,c}(L;\R)=0$, one obtains an
identification $H^1_c(L;\R)\simeq
E_0\oplus\mbox{Ker}\left(H^1_{\bullet,c}(L)\stackrel{\rho}{\rightarrow}
H^1(\Sigma_0)\right)$. Comparing these identifications yields an identification
$\widetilde{H}^1_c(L;\R)\oplus
d(E')\simeq\mbox{Ker}\left(H^1_{\bullet,c}(L)\stackrel{\rho}{\rightarrow}
H^1(\Sigma_0)\right)$. The claim follows. 

Now consider the CS/AC case with $\boldsymbol{\lambda}>0$. The long exact
sequence \ref{eq:cohomsequence} with $\Sigma=\Sigma_0$ yields
\begin{equation}\label{eq:cohomsequence_accsbis}
0\rightarrow H^0(L;\R)\rightarrow H^0(\Sigma_0;\R)\rightarrow
H^1_{c,\bullet}(L;\R)\stackrel{\gamma}{\rightarrow}
H^1(L;\R)\stackrel{\rho}{\rightarrow} H^1(\Sigma_0;\R)\rightarrow\dots
\end{equation}
This proves the final claim.
\end{proof}

\begin{remark}\label{rem:topsummary}
 Compare Equations \ref{eq:closedforms_decay_cs}, \ref{eq:closedforms_decay_ac} with the corresponding equations in the statement of Corollary \ref{cor:topsummary}. When working with AC manifolds we choose to group the two topological terms of Equation \ref{eq:closedforms_decay_ac} into one space $H^1_c(L;\R)$. When working with CS manifolds we prefer to keep the two topological terms of Equation \ref{eq:closedforms_decay_cs} separate and to emphasize the ``geometric'' meaning of one of them as kernel of a certain restriction map. These choices are based on the different roles that these spaces will play in Section \ref{s:moduli}, cf. also the ``concluding remarks'' there. 
\end{remark}


\section{Lagrangian conifolds}\label{s:lagconifolds}
A priori, an \textit{immersed conifold} or \textit{subconifold} in a Riemannian ambient space $(M,g)$ might simply be defined as an immersed
submanifold whose topology and induced metric is of the type defined in Section
\ref{s:geometry_review}. However, for the purposes of this article it is convenient
to strengthen the hypotheses by adding the requirement that the submanifold be asymptotic to a specific immersed cone at each singularity and at each AC end. If the submanifold has only conical singularities then $M$ can be any  Riemannian manifold; if the submanifold has asymptotically conical ends then, to set up the definitions, it is necessary that $M$ also have a conifold structure. 

For the sake of brevity our presentation will cover the case of immersed CS conifolds in general ambient spaces but it will discuss immersed conifolds with AC ends only in the ambient space $\C^m$, which is the ambient space of most interest to us.

We will also focus from the start on Lagrangian immersions in K\"ahler
ambient spaces, because these are the main objects of this paper.

\begin{definition} \label{def:lagr}
Let $(M^{2m},\omega)$ be a symplectic manifold. An embedded or immersed
submanifold $\iota:L^m\rightarrow M$ is \textit{Lagrangian} if
$\iota^*\omega\equiv 0$. The immersion allows us to view the tangent bundle $TL$
of $L$ as a subbundle of $TM$ (more precisely, of $\iota^*TM$). When $M$ is
K\"ahler with structures $(g,J,\omega)$ it is simple to check that $L$ is
Lagrangian iff $J$ maps $TL$ to the normal bundle $NL$ of $L$, \textit{i.e.}
$J(TL)=NL$. 
\end{definition}

\begin{definition}\label{def:aclagsub}
Let $L^m$ be a smooth manifold. Assume given a Lagrangian immersion
$\iota:L\rightarrow \C^m$, the latter endowed with its standard structures
$\tilde{J},\tilde{\omega}$. We say that $(L,\iota)$ is an \textit{asymptotically conical Lagrangian submanifold}
with \textit{rate} $\boldsymbol{\lambda}$ if it satisfies the following
conditions.
\begin{enumerate}
\item We are given a compact subset $K\subset L$ such that $S:=L\setminus K$ has
a finite number of connected components $S_1,\dots,S_e$.
\item We are given Lagrangian cones $\mathcal{C}_i\subset \C^m$ with smooth
connected links $\Sigma_i:=\mathcal{C}_i\bigcap \Sph^{2m-1}$. Let
$\iota_i:\Sigma_i\times (0,\infty)\rightarrow \C^m$ denote the natural
immersions, parametrizing $\mathcal{C}_i$.
\item We are finally given an $e$-tuple of \textit{convergence rates}
$\boldsymbol{\lambda}=(\lambda_1,\dots,\lambda_e)$ with $\lambda_i<2$, \textit{centers} $p_i\in\C^m$ and
diffeomorphisms $\phi_i:\Sigma_i\times [R,\infty)\rightarrow \overline{S_i}$
for some $R>0$ such that, for $r\rightarrow\infty$ and all $k\geq 0$,
\begin{equation}\label{eq:aclagdecay}
|\tnabla^k(\iota\circ\phi_i-(\iota_i+p_i)|=O(r^{\lambda_i-1-k})
\end{equation}
with respect to the conical metric $\tg_i$ on $\cone_i$.
\end{enumerate}
\end{definition}

Notice that the restriction $\lambda_i<2$ ensures that the cone is unique but is
weak enough to allow the submanifold to converge to a translated copy
$\mathcal{C}_i+p'_i$ of the cone (\textit{e.g.} if $\lambda_i=1$), or even to
slowly pull away from the cone (if $\lambda_i>1$). 

\begin{definition}\label{def:cslagsub}
Let $\bar{L}^m$ be a smooth manifold except for a finite number of possibly singular
points $\{x_1,\dots,x_e\}$. Assume given a continuous map
$\iota:\bar{L}\rightarrow \C^m$ which restricts to a smooth Lagrangian immersion of
$L:=\bar{L}\setminus\{x_1,\dots,x_e\}$. We say that $(\bar{L},\iota)$ or $(L,\iota)$ is a \textit{conically singular
Lagrangian submanifold} with \textit{rate} $\boldsymbol{\mu}$ if it satisfies
the following conditions.
\begin{enumerate}
\item We are given open connected neighbourhoods $S_i$ of $x_i$.
\item We are given Lagrangian cones $\mathcal{C}_i\subset \C^m$ with smooth
connected links $\Sigma_i:=\mathcal{C}_i\bigcap \Sph^{2m-1}$. Let
$\iota_i:\Sigma_i\times (0,\infty)\rightarrow \C^m$ denote the natural
immersions, parametrizing $\mathcal{C}_i$.
\item We are finally given an $e$-tuple of \textit{convergence rates}
$\boldsymbol{\mu}=(\mu_1,\dots,\mu_e)$ with $\mu_i>2$, \textit{centers} $p_i\in\C^m$ and diffeomorphisms
$\phi_i:\Sigma_i\times (0,\epsilon]\rightarrow \overline{S_i}\setminus\{x_i\}$ such that, for
$r\rightarrow 0$ and all $k\geq 0$,
\begin{equation}\label{eq:cslagdecay}
|\tnabla^k(\iota\circ\phi_i-(\iota_i+p_i))|=O(r^{\mu_i-1-k})
\end{equation}
with respect to the conical metric $\tg_i$ on $\cone_i$. Notice that our assumptions imply that $\iota(x_i)=p_i$.
\end{enumerate}
\end{definition}

It is simple to check that AC Lagrangian submanifolds, with the induced metric,
satisfy Definition \ref{def:ac_manifold} with $\nu_i=\lambda_i-2$. The analogous
fact holds for CS Lagrangian submanifolds.

\begin{definition} \label{def:accslagsub}
Let $\bar{L}^m$ be a smooth manifold except for a finite number of possibly singular
points $\{x_1,\dots,x_s\}$ and with $l$ ends. Assume
given a continuous map $\iota:\bar{L}\rightarrow \C^m$ which restricts to a smooth Lagrangian
immersion of $L:=\bar{L}\setminus\{x_1,\dots,x_s\}$. We say that $(\bar{L},\iota)$ or
$(L,\iota)$ is a \textit{CS/AC Lagrangian submanifold} with \textit{rate}
$(\boldsymbol{\mu},\boldsymbol{\lambda})$ if in a neighbourhood of the points
$x_i$ it has the structure of a CS submanifold with rates $\mu_i$ and in a
neighbourhood of the remaining ends it has the structure of an AC submanifold
with rates $\lambda_i$.

We use the generic term \textit{Lagrangian conifold} to indicate any CS, AC or CS/AC Lagrangian
submanifold.
\end{definition}

\begin{example} \label{e:lagcone}
Let $\mathcal{C}$ be a cone in $\C^m$ with smooth link $\Sigma^{m-1}$. It can be
shown that $\mathcal{C}$ is a Lagrangian iff $\Sigma$ is \textit{Legendrian} in
$\Sph^{2m-1}$ with respect to the natural \textit{contact structure} on the
sphere. Then $\mathcal{C}$ is a CS/AC Lagrangian submanifold of $\C^m$ with rate
$(\boldsymbol{\mu,\lambda})$ for any $\boldsymbol{\mu}$ and
$\boldsymbol{\lambda}$.
\end{example}

The definition of CS Lagrangian submanifolds can be generalized to K\"ahler
ambient spaces as follows. Once again we denote the standard structures on
$\C^m$ by $\tilde{J}$, $\tilde{\omega}$.
\begin{definition} \label{def:generalcslagsub}
Let $(M^{2m},J,\omega)$ be a K\"ahler manifold and $\bar{L}^m$ be a smooth
manifold except for a finite number of possibly singular points $\{x_1,\dots,x_e\}$. Assume given a continuous map $\iota:\bar{L}\rightarrow M$ which restricts
to a smooth Lagrangian immersion of $L:=\bar{L}\setminus\{x_1,\dots,x_e\}$. We say that $(\bar{L},\iota)$ or
$(L,\iota)$ is a \textit{Lagrangian submanifold with conical singularities} (CS
Lagrangian submanifold) if it satisfies the following conditions.
\begin{enumerate}
\item We are given isomorphisms $\upsilon_i:\C^m\rightarrow T_{\iota(x_i)}M$ such that
$\upsilon_i^*\omega=\tilde{\omega}$ and $\upsilon_i^*J=\tilde{J}$. 

According to Darboux' theorem, cf. \textit{e.g.} \cite{weinstein}, there then
exist an open ball $B_R$ in $\C^m$ (of small radius $R$) and diffeomorphisms
$\Upsilon_i:B_R\rightarrow M$ such that $\Upsilon(0)=\iota(x_i)$,
$d\Upsilon_i(0)=\upsilon_i$ and $\Upsilon_i^*\omega=\tilde{\omega}$.
\item We are given open neighbourhoods $S_i$ of $x_i$ in $\bar{L}$. We assume $S_i$ are small, in the sense that
the compositions
\begin{equation*}
 \Upsilon_i^{-1}\circ\iota:S_i\rightarrow B_R 
\end{equation*}
are well-defined. 

We are
also given Lagrangian cones $\mathcal{C}_i\subset\C^m$ with smooth connected
links $\Sigma_i:=\mathcal{C}_i\bigcap \Sph^{2m-1}$. Let $\iota_i:\Sigma_i\times
(0,\infty)\rightarrow \C^m$ denote the natural immersions, parametrizing
$\cone_i$. 
\item We are finally given an $e$-tuple of \textit{convergence rates}
$\boldsymbol{\mu}=(\mu_1,\dots,\mu_e)$ with $\mu_i\in (2,3)$ and diffeomorphisms 
$\phi_i:\Sigma_i\times (0,\epsilon]\rightarrow \overline{S_i}\setminus\{x_i\}$ such that, as $r\rightarrow 0$ and
for all $k\geq 0$,
\begin{equation}\label{eq:generalcslagdecay}
|\tnabla^k(\Upsilon_i^{-1}\circ\iota\circ\phi_i-\iota_i)|=O(r^{\mu_i-1-k})
\end{equation}
with respect to the conical metric $\tg_i$ on $\cone_i$.
\end{enumerate}
We call $x_i$ the \textit{singularities} of $\bar{L}$ and $\upsilon_i$ the
\textit{identifications}.
\end{definition}
One can check that, when $M=\C^m$, Definition \ref{def:generalcslagsub} coincides with Definition \ref{def:cslagsub} if we choose $\Upsilon_i(x):=x+\iota(x_i)$. Notice that the local diffeomorphisms between $M$ and $\C^m$ are prescribed only up to first order. Changing the diffeomorphism $\Upsilon_i$ (while keeping $\upsilon_i$ fixed) will perturb the map $\phi_i$ (and its derivatives) by a term of order $O(r^{2-k})$. In order to make the rate be independent of the particular diffeomorphism chosen, we need to introduce a constraint on the range of $\mu_i$ ensuring that $O(r^{2-k})<O(r^{\mu_i-1-k})$, thus $\mu_i<3$.


\subsection{Deformations of Lagrangian conifolds}\label{s:lagdefs}
We now want to understand how to parametrize the Lagrangian
deformations of a given Lagrangian conifold $L\subset M$. Since the Lagrangian
condition is invariant under reparametrization of $L$, to avoid huge amounts of
geometric redundancy it is best to work in terms of non-parametrized
submanifolds; in other words, in terms of equivalence classes of immersed
submanifolds, where two immersions are \textit{equivalent} if they differ by a
reparametrization. Then, to parametrize the possible deformations of $L$, it is
sufficient to prove a \textit{Lagrangian neighbourhood theorem}.

Recall that, given any manifold $L$, there is a \textit{tautological} 1-form $\hat{\lambda}$ on $T^*L$ defined by
$\hat{\lambda}[\alpha](v):=\alpha(\pi_*(v))$, where $\pi:T^*L\rightarrow L$ is
the natural projection. Then $\hat{\omega}:=-d\hat{\lambda}$ defines a canonical symplectic structure on $T^*L$.

The following classical result, going back to \cite{souriau} and \cite{weinstein}, is the most basic version of the Lagrangian neighbourhood theorem.
\begin{theorem}\label{th:nbd_weinstein}
Let $(M,\omega)$ be a symplectic manifold. Let $L\subset M$ be a smooth compact
Lagrangian submanifold. Then there exist a neighbourhood $\mathcal{U}$ of the
zero section of $L$ inside its cotangent bundle $T^*L$ and an embedding
$\Phi_L:\mathcal{U}\rightarrow M$ such that $\Phi_{L|L}=Id:L\rightarrow L$ and
$\Phi_L^*\omega=\hat{\omega}$. 
\end{theorem}

\begin{remark} Although the statement is for embedded submanifolds, it
is not difficult to extend it to immersed compact Lagrangian submanifolds by
working locally. In this case $\Phi_L$ will only be a local embedding. 
\end{remark}

Let $C^\infty(\mathcal{U})$ denote the space of smooth 1-forms on $L$ whose
graph lies in $\mathcal{U}$. In particular $\Phi_L$ defines by composition an injective map 
\begin{equation}\label{eq:cptlagrcorrespondence}
\Phi_L:C^\infty(\mathcal{U})\rightarrow \mbox{Imm}(L,M)/\mbox{Diff}(L).
\end{equation}
An important point about this map is that any submanifold which admits a
parametrization which is $C^1$-close to some parametrization of $L$ belongs to
the image of $\Phi_L$, \textit{i.e.} corresponds to a 1-form $\alpha$. 

One can check that a section $\alpha\in C^\infty(\mathcal{U})$ is closed iff the corresponding submanifold $\Phi_L\circ\alpha$ is Lagrangian. This allows us to specialize the correspondence of Equation \ref{eq:cptlagrcorrespondence} to Lagrangian immersions. In particular, let $\mbox{Lag}(L,M)$ denote the set of Lagrangian immersions from $L$ into $M$.
Using the Fr\'echet topology on
$C^\infty(\mathcal{U})$ one can locally define a topology on
$\mbox{Lag}(L,M)/\mbox{Diff}(L)$; on the intersection of
any two open sets these topologies coincide, so we obtain a global topology on
$\mbox{Lag}(L,M)/\mbox{Diff}(L)$. The connected component containing the given
$L\subset M$ defines the \textit{moduli space of Lagrangian deformations of
$L$}. Coupling Theorem \ref{th:nbd_weinstein} with Decomposition
\ref{decomp:closedforms} of Section \ref{ss:closedforms} gives a good idea of the local structure of this space.

In \cite{joyce:I}, Joyce set up an analogous framework for dealing with deformations of Lagrangian conifolds. In this case it is necessary to also control the rates of convergence of the deformations, using the rates of convergence of the closed 1-forms. This requires a very careful choice of symplectomorphism $\Phi_L$ along the ends of $L$. The reader can find a detailed explanation of how to do this in \cite{pacini:sldefs}. The final result is as follows. 

\begin{theorem}\label{th:lagr_defs}
Let $L\subset \C^m$ be a
Lagrangian conifold in $\C^m$ with asymptotic cone $\mathcal{C}$ and rate $(\boldsymbol{\mu},\boldsymbol{\lambda})$. Then there exist a neighbourhood $\mathcal{U}$ of the
zero section of $L$ inside its cotangent bundle $T^*L$ and an embedding
$\Phi_L:\mathcal{U}\rightarrow \C^m$ such that $\Phi_{L|L}=Id:L\rightarrow L$ and 
$\Phi_L^*\omega=\hat{\omega}$.

For any weight $\boldsymbol{\beta}$, let $C^\infty_{\boldsymbol{\beta}}(\mathcal{U})$ denote the corresponding space of smooth 1-forms on $L$ whose
graph lies in $\mathcal{U}$.
A section $\alpha\in C^\infty_{(\boldsymbol{\mu}-1,\boldsymbol{\lambda}-1)}(\mathcal{U})$ is closed if and only if the corresponding immersion $\Phi_L\circ\alpha$ is a Lagrangian conifold with the same asymptotic cone $\mathcal{C}$ and rate $(\boldsymbol{\mu},\boldsymbol{\lambda})$. 
\end{theorem}

In complete analogy with the compact case, we can use Theorem \ref{th:lagr_defs} to define a topology on the set of
Lagrangian conifolds which admit a parametrization
$\iota:L\rightarrow\C^m$ which is asymptotic to $\mathcal{C}$
with rate $(\boldsymbol{\mu},\boldsymbol{\lambda})$. The connected component containing a given $\iota$
defines the \textit{moduli space of Lagrangian deformations of
$(L,\iota)$ with rate $(\boldsymbol{\mu},\boldsymbol{\lambda})$}. 

Coupling these results with Decompositions \ref{decomp:closedforms_growth},
\ref{decomp:closedforms_decay} and \ref{decomp:closedforms_accs} of Section \ref{ss:closedforms_bis} now gives a
good idea of the local structure of the corresponding moduli spaces of
Lagrangian deformations of $\iota$.

Up to here, the given Lagrangian conifold $(L,\iota)$ has been deformed
keeping the singular points $\{\iota(x_1),\dots,\iota(x_s)\}$ fixed in the ambient manifold $\C^m$. It is
also natural to want to deform $L$ allowing the singular points to move in $\C^m$. Analogously, one might want to allow the corresponding
Lagrangian cones $\cone_i$ to rotate. The correct set-up for doing
this is as follows. The ideas are based on \cite{joyce:II}
Section 5.1. Define
\begin{equation}\label{eq:P}
P:=\{(p,\upsilon):p\in \C^m,\ \upsilon\in \unitary m\}.
\end{equation}
$P$ is a $\unitary m$-principal fibre bundle over $\C^m$ with the action 
$$\unitary m\times P\rightarrow P,\ \ M\cdot (p,\upsilon):=(p,\upsilon\circ
M^{-1}).$$
As such, $P$ is a smooth manifold of dimension $m^2+2m$.

Our aim is to use one copy of $P$ to parametrize the location of each singular
point $p_i=\iota(x_i)\in \C^m$ and the corresponding asymptotic cone $\upsilon_i(\mathcal{C}_i)$: the group action will allow the cone to rotate leaving the singular point
fixed. As we are interested only in small deformations of $L$ we can restrict
our attention to a small open neighbourhood of the pair $(p_i,Id)\in P$. In general the
$\cone_i$ will have some symmetry group $G_i\subset \unitary m$, \textit{i.e.}
the action of this $G_i$ will leave the cone fixed. To ensure that we have no
redundant parameters we must therefore further restrict our attention to a
\textit{slice} of our open neighbourhood, \textit{i.e.} a smooth submanifold
transverse to the orbits of $G_i$. We denote this slice $\E_i$: it is a subset
of $P$ containing $(p_i,Id)$ and of dimension $m^2+2m-\mbox{dim}(G_i)$.
We then set $\E:=\E_1\times\dots\times\E_s$. The point
$e:=(p_1,Id),\dots,(p_s,Id))\in \E$ corresponds to the initial data.

We now want to extend the initial datum of $(L,\iota)$ to a family of Lagrangian submanifolds
$(L,\iota_{\tilde{e}})$ parametrized by
$\tilde{e}=((\tilde{p}_1,\tilde{\upsilon}_1),\dots,(\tilde{p}_s,\tilde{\upsilon}_s))\in
\E$. Each $(L,\iota_{\tilde{e}})$ should satisfy $\iota_{\tilde{e}}(x_i)=\tilde{p_i}$ and have  asymptotic cones $\tilde{\upsilon}_i(\mathcal{C}_i)$. We further require that $\iota_e=\iota$ globally and that $\iota_{\tilde{e}}=\iota$ outside a neighbourhood of the singularities.
The construction of such a family is actually straight-forward: for each $\tilde{e}$, it reduces to a choice of a
compactly-supported symplectomorphism of $\C^m$ (which we continue to denote $\tilde{e}$) which, near each $p_i$, extends the maps $x\mapsto \tilde{p}_i+\tilde{\upsilon}_i(x-p_i)$. We then obtain immersions $\iota_{\tilde{e}}:=\tilde{e}\circ\iota$ and 
embeddings $\Phi_{L}^{\tilde{e}}:=\tilde{e}\circ\Phi_L:\mathcal{U}\rightarrow \C^m$ which, away from the
singularities, coincide with $\iota$ and $\Phi_{L}$. The final result is that, after such a choice, the
\textit{moduli space of Lagrangian deformations of $L$ with rate
$(\boldsymbol{\mu},\boldsymbol{\lambda})$ and moving singularities} can be parametrized in terms of
pairs $(\tilde{e}, \alpha)$ where $\tilde{e}\in\E$ and $\alpha$ is a closed 1-form
on $L$ belonging to the space $C^\infty_{(\boldsymbol{\mu}-1,\boldsymbol{\lambda}-1)}(\mathcal{U})$. 

\begin{remark}\label{rem:Mextension}
All the above results and constructions can be extended to CS submanifolds in $M$, using appropriate compositions
by $\Upsilon_i$. In
this case we set $P:=\{(p,\upsilon)\}$, where $p\in M$ and
$\upsilon: \C^m\rightarrow T_pM$ such that $\upsilon^*\omega=\tilde{\omega}$, $\upsilon^*J=\tilde{J}$.
\end{remark}


\section{Special Lagrangian conifolds}\label{s:slgeometry}
\begin{definition}\label{def:cy} A \textit{Calabi-Yau} (CY) manifold is the data
of a K\"ahler manifold ($M^{2m}$,$g$,$J$,$\omega$) and a non-zero ($m,0$)-form
$\Omega$ satisfying $\nabla\Omega\equiv 0$ and normalized by the condition
$\omega^m/m!=(-1)^{m(m-1)/2}(i/2)^m\Omega\wedge\bar{\Omega}$.

In particular $\Omega$ is holomorphic and the holonomy of $(M,g)$ is contained
in $\sunitary m$. We will refer to $\Omega$ as the \textit{holomorphic volume
form} on $M$.
\end{definition}

\begin{definition}\label{def:sl}
Let $M^{2m}$ be a CY manifold and $L^m\rightarrow M$ be an immersed or embedded
Lagrangian submanifold. We can restrict $\Omega$ to $L$, obtaining a
non-vanishing complex-valued $m$-form $\Omega_{|L}$ on $L$. We say that $L$ is
\textit{special Lagrangian} (SL) iff this form is real, \textit{i.e.}
$\Imag\,\Omega_{|L}\equiv 0$. In this case $\Real\,\Omega_{|L}$ defines a volume
form on $L$, thus a natural orientation.
\end{definition}

\begin{example}\label{e:C^m}
The simplest example of a CY manifold is $\C^m$ with its standard structures
$\tilde{g}$, $\tilde{J}$, $\tilde{\omega}$ and
$\tilde{\Omega}:=dz^1\wedge\dots\wedge dz^m$. The simplest example of SL submanifold in $\C^m$ is the standard plane $\R^m$. Any other SL plane $\Pi$ can be obtained by rotating $\R^m$ via a matrix in $\sunitary m$. Using this fact, it is simple to show that, for any normal vector $v\in\Pi^\perp$,
$(i_v\Imag\,\tilde{\Omega})_{|\Pi}=-\star\alpha$, where $\alpha=\iota_v\tilde{\omega}=\tilde{\omega}(v,\cdot)_{|\Pi}$ and $\star$ is the Hodge star operator. 

Now let $L$ be a general SL submanifold in a general CY manifold $M$. Fixing a point $x\in L$, one can choose an isomorphism
$T_xM\simeq \C^m$ identifying the CY structures on $T_xM$ with the standard
structures on $\C^m$. This map will identify $T_xL$ with a SL $m$-plane $\Pi$ in
$\C^m$, showing that the above relationship holds pointwise for $L$. The final result is the useful formula
\begin{equation}\label{eq:twist}
(i_v\Imag\,\Omega)_{|L}=-\star\alpha,
\end{equation}
for any normal vector $v\in TL^\perp$ and $\alpha=\iota_v\omega_{|TL}$.
\end{example}

\begin{definition} \label{def:accssl}
We can define AC, CS and CS/AC special Lagrangian submanifolds in $\C^m$ exactly
as in Definitions \ref{def:aclagsub}, \ref{def:cslagsub} and
\ref{def:accslagsub}, simply adding the requirement that the submanifolds be
special Lagrangian. In particular this implies that the cones $\mathcal{C}_i$
are SL in $\C^m$. Following Definition \ref{def:generalcslagsub} we can also
define CS special Lagrangian submanifolds in a general CY manifold $M$: in this
case it is necessary to also add the requirement that
$\upsilon_i^*\Omega=\tilde{\Omega}$.

We use the generic term \textit{special Lagrangian conifold} to refer to any of
the above.
\end{definition}

\begin{remark} \label{rem:muregularity} 
It follows from Joyce \cite{joyce:I} Theorem 5.5 that if $L$ is a CS or CS/AC SL
submanifold with respect to some rate $\boldsymbol{\mu}=2+\epsilon$ with
$\epsilon$ in a certain range $(0,\epsilon_0)$ then it is also CS or CS/AC with
respect to any other rate of the form $\boldsymbol{\mu}'=2+\epsilon'$ with
$\epsilon'\in (0,\epsilon_0)$. The precise value of $\epsilon_0$ is determined
by the exceptional weights of the cones $\mathcal{C}_i$, as in Section \ref{s:reviewlaplace}. We refer to \cite{joyce:I} for
details.
\end{remark}
  
SL submanifolds are \textit{calibrated submanifolds} in the sense of \cite{harveylawson}. This implies that they are volume-minimizing in their homology class, and in particular are minimal. It is well-known that the ambient space $\C^m$ cannot contain compact minimal submanifolds. It follows that any SL conifold in $\C^m$ must have at least one AC end. In other words, there is no point in studying CS SLs in $\C^m$.

\begin{example} \label{e:slcone}
Let $\mathcal{C}$ be a Lagrangian cone in $\C^m$ with smooth link
$\Sigma^{m-1}$. It can be shown that $\mathcal{C}$ is SL (with respect to some
holomorphic volume form $e^{i\theta}\tilde{\Omega}$) iff $\Sigma$ is minimal in
$\Sph^{2m-1}$ with respect to the natural metric on the sphere, so 
$\mathcal{C}$ is a CS/AC SL in $\C^m$. We refer to \cite{harveylawson},
\cite{haskins}, \cite{haskinskapouleas}, \cite{haskinspacini}, \cite{joyce:symmetries} for examples.

We refer to Joyce \cite{joyce:V} Section 6.4 for examples of AC SLs in $\C^m$
with various rates.
\end{example}

Lagrangian submanifolds (especially the immersed ones) tend to be very ``soft"
objects: for example, Section \ref{s:lagdefs} shows that they have
infinite-dimensional moduli spaces. They also easily allow for cutting, pasting
and desingularization procedures. The ``special" condition rigidifies them
considerably: the corresponding gluing and desingularization
processes require much ``harder" techniques. We refer to 
\cite{haskinskapouleas}, \cite{joyce:III}, \cite{joyce:IV}, \cite{pacini:slgluing} for recent gluing
results and \cite{haskinspacini} for local desingularization issues. The main goal of this paper is to ``quantify'' this notion of rigidity by examining the problem of SL deformations and calculating the corresponding degrees of freedom.


\subsection{Setting up the SL deformation problem}\label{s:setup}
If $\iota:L\rightarrow M$ is a SL submanifold we can specialize the framework of Section
\ref{s:lagdefs} to study the SL deformations of $L$. Notice that the SL
condition is again invariant under reparametrizations. Thus, if $L$ is smooth
and compact, the \textit{moduli space} $\mathcal{M}_L$ \textit{of SL
deformations of $(L,\iota)$} can be defined as the connected component containing $\iota$ of
the subset of SL immersions in $\mbox{Lag}(L,M)/\mbox{Diff}(L)$. As in Section \ref{s:lagdefs}, if $L$ is a SL conifold with specific rates of growth/decay on the ends we can obtain
moduli spaces of SL deformations of $(L,\iota)$ with those same rates by
restricting our attention to closed 1-forms on $L$ which satisfy
corresponding growth/decay conditions. 

Our ultimate goal is to determine situations in which moduli spaces of SL conifolds admit a
natural smooth structure with respect to which they are finite-dimensional
manifolds. In particular, we need to identify the obstructions which may prevent this
from happening. Generally speaking, the strategy for proving these results will
be to view $\mathcal{M}_L$ locally as the zero set of some smooth map $F$
defined on the space of closed forms in $C^\infty(\mathcal{U})$ (when $L$ is
smooth and compact) or in
$C^\infty_{(\boldsymbol{\mu}-1,\boldsymbol{\lambda}-1)}(\mathcal{U})$ (when $L$
is CS/AC with rate $(\boldsymbol{\mu},\boldsymbol{\lambda})$): we can then
attempt to use the Implicit Function Theorem to prove that this zero set
is smooth. 

The choice of $F$ is dictated by Definition \ref{def:sl}: basically, if $\Omega$ denotes
the given holomorphic volume form on $M$ then $F$ must compute the values
of $\Imag\,\Omega$ on each Lagrangian deformation of $L$. 

\ 

\textit{Note: }To simplify the notation, from now on we will drop the immersion $\iota:L\rightarrow M$ and simply identify $L$ with its image. In particular we will identify the singularities $x_i$ with their images $\iota(x_i)$. 

\ 

As a first case, let $L\subset M$ be a smooth compact SL submanifold, endowed with the induced
metric $g$ and orientation. Let $\star$ denote the Hodge
star operator defined on $L$ by $g$ and the orientation. Define $\Phi_L:\mathcal{U}\rightarrow M$ as in
Section \ref{s:lagdefs}. Let $\mathcal{D}_L$ denote the space of closed 1-forms on $L$ whose graph lies
in $\mathcal{U}$. We then define the map $F$ as follows. 
\begin{equation}\label{eq:defF}
F:\mathcal{D}_L\rightarrow C^\infty (L),\ \ \alpha\mapsto
\star(\alpha^*(\Phi_L^*\Imag\,\Omega))=\star((\Phi_L\circ\alpha)^*\Imag\,\Omega).
\end{equation}

The following result is due to \cite{mclean}.
\begin{prop} \label{prop:cptnonlinear}
The non-linear map $F$ has the following properties:
\begin{enumerate}
\item The set $F^{-1}(0)$ parametrizes the space of all SL deformations of $L$
which are $C^1$-close to $L$.
\item $F$ is a smooth map between Fr\'echet spaces.  Furthermore, for each
$\alpha\in\mathcal{D}_L$, $\int_LF(\alpha)\,vol_g=0$.
\item The linearization $dF[0]$ of $F$ at $0$ coincides with the operator $d^*$,
\textit{i.e.} 
\begin{equation}\label{eq:linearization}
dF[0](\alpha)=d^*\alpha.
\end{equation}
\end{enumerate}
\end{prop}
\begin{proof} It is instructive to sketch the proof of Equation \ref{eq:linearization}. We refer to \cite{pacini:sldefs} for full  details. To simplify the notation, we identify $\mathcal{U}$ with its image in $M$ via $\Phi_L$.

Fix any $\alpha\in \Lambda^1(L)$. The Lagrangian condition implies that the vector field $v$ defined along $L$ by imposing
$\alpha(\cdot)\equiv\omega(v,\cdot)$ is normal to $TL$. We can extend $v$ to a global vector field
$v$ on $M$. Let $\phi_s$ denote any 1-parameter family of diffeomorphisms of $M$
such that $d/ds(\phi_s(x))_{|s=0}=v(x)$. Then the two 1-parameter families of
$m$-forms on $L$, $(s\alpha)^*(\Imag\,\Omega)=\pi_*(\Imag \,\Omega_{|\Gamma(s\alpha)})$ and
$(\phi_s^*\Imag\,\Omega)_{|L}$, coincide up to first order so that standard
calculus of Lie derivatives shows that
\begin{eqnarray*}
dF[0](\alpha)\,vol_g &=& d/ds
(F(s\alpha)\,vol_g)_{|s=0}\\
&=& d/ds(\phi_s^*\Imag\,\Omega)_{|L;\,s=0}\\
&=& (\mathcal{L}_v\Imag\,\Omega)_{|L}= (di_v\Imag\,\Omega)_{|L},
\end{eqnarray*}
where in the last equality we use \textit{Cartan's formula}
$\mathcal{L}_v=di_v+i_vd$ and the fact that $\Imag\,\Omega$ is closed. We now apply Equation \ref{eq:twist} to conclude.
\end{proof}

Our main goal is to understand how to parametrize the SL deformations of a SL conifold $L$ in $\C^m$. As in Section \ref{s:lagdefs}, we want to allow the singularities of $L$ to move. The SL constraint suggests that we modify the definition given in Equation \ref{eq:P} as follows:
\begin{equation}
\tilde{P}:=\{(p,\upsilon):p\in \C^m,\ \upsilon\in \sunitary m\}.
\end{equation}
$\tilde{P}$ is then a $\sunitary m$-principal fibre bundle over $\C^m$ of
dimension $m^2+2m-1$. For each end, the cone $\cone_i$ will have symmetry
group $G_i\subset \sunitary m$. Let
$\tilde{\E}_i$ denote a smooth submanifold of $\tilde{P}$ transverse to the
orbits of $G_i$. It has dimension $m^2+2m-1-\mbox{dim}(G_i)$. Set
$\tilde{\E}:=\tilde{\E}_1\times\dots\times\tilde{\E}_s$. We then define
Lagrangian conifolds $(L, \iota_{\tilde{e}})$ and embeddings $\Phi_{L}^{\tilde{e}}$ as before.

Let $\mathcal{D}_{L}$
denote the space of closed 1-forms in
$C^\infty_{(\boldsymbol{\mu}-1,\boldsymbol{\lambda}-1)}(\Lambda^1)$ whose graph
lies in $\mathcal{U}$. Consider the map
\begin{equation}\label{eq:csacdefF}
F:\tilde{\E}\times\mathcal{D}_{L}\rightarrow C^\infty_{(\boldsymbol{\mu}-2,\boldsymbol{\lambda}-2)} (L),\
\ (\tilde{e},\alpha)\mapsto
\star(\alpha^*(\Phi_{L}^{\tilde{e}*}\Imag\,\Omega)).
\end{equation}

\begin{prop}\label{prop:accsnonlinear}
Let $L$ be a SL conifold in $\C^m$. Then the map $F$ has the following properties:
\begin{enumerate}
\item The set $F^{-1}(0)$ parametrizes the space of all SL deformations of $L$
which are $C^1$-close to $L$ away from the singularities, have centers $\tilde{p}_i$ and are asymptotic to
$\tilde{\upsilon}_i(\mathcal{C}_i)$ with rate $(\boldsymbol{\mu},\boldsymbol{\lambda})$ for some
choice of $(\tilde{p}_i,\tilde{\upsilon}_i)$ near $(p_i,\upsilon_i)$.
\item $F$ is a (locally) well-defined smooth map between Fr\'echet spaces. In
particular, for each $\alpha\in\mathcal{D}_{L}$, $F(\alpha)\in
C^\infty_{(\boldsymbol{\mu}-2,\boldsymbol{\lambda}-2)}(L)$.
\item There exists an injective linear map $\chi:T_e\tilde{\E}\rightarrow
C^\infty_{\boldsymbol{0}}(L)$ such that (i) $\chi(y)\equiv 0$ away from the
singularities and (ii) the linearized map $dF[0]:T_e\tilde{\E}\oplus
C^\infty_{(\boldsymbol{\mu}-1,\boldsymbol{\lambda}-1)}(\Lambda^1)\rightarrow
C^\infty_{(\boldsymbol{\mu}-2,\boldsymbol{\lambda}-2)}(L)$ satisfies
\begin{equation}\label{eq:linearization_bis}
dF[0](y,\alpha)=\Delta_g\,\chi(y)+d^*\alpha.
\end{equation}
\end{enumerate}
\end{prop}

\begin{proof} As for Equation \ref{eq:linearization}, it is instructive to at least sketch the proof of Equation \ref{eq:linearization_bis}. Again we refer to \cite{pacini:sldefs} for full details.

The linearization of $F$ with respect to
directions in $C^\infty_{(\boldsymbol{\mu}-1,\boldsymbol{\lambda}-1)}(\Lambda^1)$ can be computed as in
Proposition \ref{prop:cptnonlinear}. Now choose $y\in T_e\tilde{\E}$
corresponding to a curve $\tilde{e}_s\in\tilde{\E}$ such that $\tilde{e}_0=e$. Recall from the paragraph immediately preceding Remark \ref{rem:Mextension} that we can identify $\tilde{e}_s$ with a curve of compactly-supported symplectomorphisms of $\C^m$ which, near each singularity, extend the action of $\sunitary m\ltimes\C^m$ on $\C^m$. The tangent direction $y$ can then be identified with the vector field induced by $\tilde{e}_s$ on $\C^m$, \textit{i.e.} $y=d/ds(\tilde{e}_s)_{|s=0}$. Then, as in
Proposition \ref{prop:cptnonlinear} and with the same identifications, 
\begin{eqnarray}\label{eq:linearization_ter}
dF[0](y)\,vol_g &=&
d/ds(F(\tilde{e}_s,0)\,vol_g)_{|s=0}=d/ds((\tilde{e}_s)^*\Imag\,\tilde{\Omega})_{|L;s=0}\\
&=& (\mathcal{L}_y\Imag\,\tilde{\Omega})_{L}=(di_y\Imag\,\tilde{\Omega})_{|L}\nonumber\\
&=&-d\star\alpha,\nonumber
\end{eqnarray}
where $\alpha:=\iota_y\tilde{\omega}_{|L}$ is a closed 1-form on $L$. 

We now want to look more closely at this 1-form $\alpha$ near the singularities of $L$, where $\tilde{e}_s$ is a 1-parameter curve in the
group $\sunitary m\ltimes\C^m$. The action of $\sunitary m\ltimes\C^m$ on $\C^m$ admits a
\textit{moment map} $\mu:\C^m\rightarrow (Lie(\sunitary m\ltimes\C^m))^*$.
Recall that this means that $\mu$ is equivariant and that, for all $w\in
Lie(\sunitary m\ltimes\C^m)$, the corresponding function
$\mu_w:\C^m\rightarrow\R$ satisfies $d\mu_w=i_w\tilde{\omega}$, \textit{i.e.}
$w$ is a Hamiltonian vector field with Hamiltonian function $\mu_w$. The moment
map can be written explicitly, cf. \textit{e.g.} \cite{haskinspacini} Section
2.6, showing that each $\mu_w$ is at most a quadratic polynomial on $\C^m$.
Since our vector field $y$ is, near each singularity, an element of $Lie(\sunitary
m\ltimes\C^m)$ we can set $\chi(y):=\mu_y$ so that $\alpha=d\,\chi(y)$. This shows in particular that $\alpha$ is exact on the ends of $L$. Since the symplectomorphisms $\tilde{e}_s$ have compact support away from the singularities, we see that $\alpha\equiv 0$ on $K\subset L$. The long exact sequence \ref{eq:cohomsequence} then shows that $\alpha$ is globally exact so we can write $\alpha=d\,\chi(y)$, for some extension $\chi(y):L\rightarrow\R$. Plugging this into Equation \ref{eq:linearization_ter} proves Equation \ref{eq:linearization_bis}. Our explicit description of $\chi(y)$ on the ends shows that it is bounded as $r\rightarrow
0$ and has lowest order terms of order 0 so $\chi(y)\in C^0_{\boldsymbol{0}}(L)$. Further calculations
show that $\chi(y)\in C^\infty_{\boldsymbol{0}}(L)$, as claimed.

For future reference we add that, for any SL submanifold $L$ in $\C^m$, Equation \ref{eq:twist} shows that
\begin{equation*}
\Delta_g(\mu_{w|L})=d^*(d\mu_{w|L})=-\star
d\star(i_w\tilde{\omega}_{|L})=\star (d
i_w\Imag\,\tilde{\Omega})_{|L}=\star(\mathcal{L}_w\Imag\,\tilde{\Omega
})_{|L}=0,
\end{equation*}
\textit{i.e.} each $\mu_w$ restricts to a harmonic function on each SL submanifold. In particular this calculation shows that $\Delta_g\,\chi(y)$ vanishes near each singularity. 
\end{proof}

If the spaces $C^\infty(L)$, $C^\infty_{\boldsymbol{\beta}}(L)$ were Banach
spaces and the relevant maps were Fredholm, we could now apply the Implicit
Function Theorem to conclude that the sets $F^{-1}(0)$, and thus
$\mathcal{M}_L$, are smooth. As however they are actually only Fr\'echet spaces,
it is instead necessary to first take the Sobolev space completions of these
spaces, then study the Fredholm properties of the linearized maps. We do this in
Section \ref{s:moduli}.

\subsection{Stable SL cones} Given  a SL conifold $L$ we will see that smoothness of $\mathcal{M}_L$ requires
an additional ``stability" assumption on the asymptotic SL cones corresponding to the conical singularities. Roughly speaking, it
is required that these cones admit no additional harmonic
functions with prescribed growth, beyond those which necessarily exist for the
geometric reasons described in the proof of Proposition \ref{prop:accsnonlinear}. No condition will be required on the asymptotic SL cones corresponding to the AC ends.

\begin{definition}\label{def:stability}
Let $\mathcal{C}$ be a SL cone in $\C^m$. Let $(\Sigma,g')$ denote the link of
$\mathcal{C}$ with the induced metric. Assume $\mathcal{C}$ has a unique
singularity at the origin; equivalently, assume that $\Sigma$ is smooth and that
it is not a sphere $\Sph^{m-1}\subset \Sph^{2m-1}$. Recall from the proof of
Proposition \ref{prop:accsnonlinear} that the standard action of $\sunitary
m\ltimes\C^m$ on $\C^m$ admits a moment map $\mu$ and that the components of
$\mu$ restrict to harmonic functions on $\mathcal{C}$. Let $G$ denote the
subgroup of $\sunitary m$ which preserves $\mathcal{C}$. Then $\mu$ defines on
$\mathcal{C}$ $2m$ linearly independent harmonic functions of linear growth; in
the notation of Definition \ref{def:exceptionalweights} these functions are
contained in the space $V_\gamma$ with $\gamma=1$. The moment map also defines
on $\mathcal{C}$ $m^2-1-\mbox{dim}(G)$ linearly independent harmonic functions
of quadratic growth: these belong to the space $V_\gamma$ with $\gamma=2$.
Constant functions define a third space of homogeneous harmonic functions on
$\mathcal{C}$, \textit{i.e.} elements in $V_\gamma$ with $\gamma=0$. In
particular, these three values of $\gamma$ are always exceptional values for the
operator $\Delta_{\tilde{g}}$ on any SL cone, in the sense of Definition
\ref{def:exceptionalweights}.

We say that $\mathcal{C}$ is \textit{stable} if these are the only functions in
$V_\gamma$ for $\gamma=0,1,2$ and if there are no other exceptional values
$\gamma$ in the interval $[0,2]$. More generally, let $L$ be a CS or CS/AC SL
submanifold. We say that a singularity $x_i$ of $L$ is \textit{stable} if the
corresponding cone $\mathcal{C}_i$ is stable.
\end{definition}

Stability is a strong condition and very few examples of stable SL cones are known. We refer to \cite{haskins:complexity}, \cite{joyce:II} and \cite{ohnita} for more details and examples.


\section{Moduli spaces of special Lagrangian conifolds}\label{s:moduli}

Recall the statement of the Implicit Function Theorem.
\begin{theorem} \label{thm:IFT}
Let $F:E_1\rightarrow E_2$ be a smooth map between Banach spaces such that
$F(0)=0$. Assume $P:=dF[0]$ is surjective and  $\mbox{Ker}(P)$ admits a closed
complement $Z$, \textit{i.e.} $E_1=\mbox{Ker}(P)\oplus Z$. Then there exists a
smooth map $\Phi:\mbox{Ker}(P)\rightarrow Z$ such that $F^{-1}(0)$ coincides
locally with the graph $\Gamma(\Phi)$ of $\Phi$. In particular, $F^{-1}(0)$ is
(locally) a smooth Banach submanifold of $E_1$.
\end{theorem}

The following result is straight-forward. 

\begin{prop}\label{prop:fredholmreducestofinitedim}
Let $F:E_1\rightarrow E_2$ be a smooth map between Banach spaces such that
$F(0)=0$. Assume $P:=dF[0]$ is Fredholm. Set $\mathcal{I}:=\mbox{Ker}(P)$ and
choose $Z$ such that $E_1=\mathcal{I}\oplus Z$. Let $\mathcal{O}$ denote a
finite-dimensional subspace of $E_2$ such that $E_2=\mathcal{O}\oplus
\mbox{Im}(P)$. Define
\begin{equation*}
G:\mathcal{O}\oplus E_1\rightarrow E_2,\ \ (\gamma,e)\mapsto \gamma+F(e).
\end{equation*}
Identify $E_1$ with $(0,E_1)\subset \mathcal{O}\oplus E_1$. Then:
\begin{enumerate}
\item The map $dG[0]=Id\oplus P$ is surjective and
$\mbox{Ker}(dG[0])=\mbox{Ker}(P)$. Thus, by the Implicit Function Theorem, there
exist $\Phi:\mathcal{I}\rightarrow\mathcal{O}\oplus Z$ such that
$G^{-1}(0)=\Gamma(\Phi)$.
\item $F^{-1}(0)=\{(i,\Phi(i)):\Phi(i)\in
Z\}=\{(i,\Phi(i)):\pi_{\mathcal{O}}\circ\Phi(i)=0\}$,
where $\pi_{\mathcal{O}}:\mathcal{O}\oplus Z\rightarrow\mathcal{O}$ is the
standard projection.
\item Let $\pi_{\mathcal{I}}:\mathcal{I}\oplus Z\rightarrow \mathcal{I}$ denote
the standard projection. Then $\pi_{\mathcal{I}}$ is a continuous open map so it
restricts to a homeomorphism
\begin{equation*}
\pi_{\mathcal{I}}:F^{-1}(0)\rightarrow (\pi_{\mathcal{O}}\circ\Phi)^{-1}(0)
\end{equation*}
between $F^{-1}(0)$ and the zero set of the smooth map
$\pi_{\mathcal{O}}\circ\Phi:\mathcal{I}\rightarrow\mathcal{O}$, which is defined
between finite-dimensional spaces.
\end{enumerate}
\end{prop}

We now have all the ingredients necessary to study the smoothness of the SL moduli space of a given SL conifold $L$. 
Equation
\ref{eq:csacdefF} described this moduli space as the zero set of a map $F$. To be able to apply the Implicit Function Theorem it is necessary to reformulate this description using Banach spaces. To this end, choose $k\geq 3$ and $p>m$ so that
$W^p_{k-1,(\boldsymbol{\mu}-1,\boldsymbol{\lambda}-1)}(\Lambda^1)\subset
C^1_{(\boldsymbol{\mu}-1,\boldsymbol{\lambda}-1)}(\Lambda^1)$. Let $\mathcal{D}_L$ denote the space
of closed 1-forms in $W^p_{k-1,(\boldsymbol{\mu}-1,\boldsymbol{\lambda}-1)}(\Lambda^1)$ whose graph
$\Gamma(\alpha)$ lies in $\mathcal{U}$. Consider the map 
\begin{equation}\label{eq:csacFsobolev}
F:\tilde{\E}\times\mathcal{D}_L\rightarrow W^p_{k-2,(\boldsymbol{\mu}-2,\boldsymbol{\lambda}-2)}(L), \ \
(\tilde{e},\alpha)\mapsto
\star(\alpha^*(\Phi_{L}^{\tilde{e}*}\Imag\,\Omega)).
\end{equation}
Since $\boldsymbol{\mu}>2$ and $\boldsymbol{\lambda}<2$, Theorem
\ref{thm:embedding} shows that
$W^p_{k-2,(\boldsymbol{\mu}-2,\boldsymbol{\lambda}-2)}(L)$ is closed under
multiplication. As in Proposition
\ref{prop:accsnonlinear}, this shows that $F$ is a (locally well-defined) smooth
map between Banach spaces with differential $dF[0](y,\alpha)=\Delta_g\,\chi(y)+d^*\alpha$. Assume $F(\tilde{\epsilon},\alpha)=0$. Regularity results of Joyce \cite{joyce:I} can then be used to show that $\alpha\in C^\infty_{(\boldsymbol{\mu}-1,\boldsymbol{\lambda}-1)}(\Lambda^1)$ so 
$F^{-1}(0)$ is locally homeomorphic, via $\Phi_L$, to $\mathcal{M}_L$. 

Notice that $F$ is a first-order map acting (up to the finite-dimensional space $\tilde{\E}$) on 1-forms. 
To prove our result it actually is useful to modify the map $F$ one more time, emphasizing the subspace of exact 1-forms: this can be achieved by switching to a second-order map acting on functions. In the course of the proof we will thus define a new map of the form
\begin{equation}\label{eq:abstract_tildeF}
\tilde{F}:K\times
W^p_{k,(\boldsymbol{\mu},\boldsymbol{\lambda})}(L)\rightarrow W^p_{k-2,(\boldsymbol{\mu}-2,\boldsymbol{\lambda}-2)}(L),
\end{equation}
 where $K$ is a
finite-dimensional vector space defined in terms of spaces introduced in
Sections \ref{ss:closedforms} and \ref{s:setup}. Geometrically, this new point of view 
corresponds to separating the obvious Hamiltonian deformations of $L$ from a
finite-dimensional space of other Lagrangian deformations. This has two benefits: (i) it allows us to make full use of the (relatively simple) theory of the Laplace operator on functions, and (ii) it emphasizes the different role played by each space.

\begin{theorem} \label{thm:accssl}
Let $L$ be a SL conifold in $\C^m$ with $s$ CS ends, $l$ AC ends
and rate $(\boldsymbol{\mu},\boldsymbol{\lambda})$. Let $\mathcal{M}_L$ denote
the moduli space of SL deformations of $L$ with moving singularities and rate
$(\boldsymbol{\mu},\boldsymbol{\lambda})$. Assume
$(\boldsymbol{\mu},\boldsymbol{\lambda})$ is non-exceptional for the Laplace operator
\begin{equation}\label{eq:accslaplacian}
\Delta_{\boldsymbol{\mu},\boldsymbol{\lambda}}:W^p_{k,(\boldsymbol{\mu},\boldsymbol{\lambda})}(L)\rightarrow
W^p_{k-2,(\boldsymbol{\mu}-2,\boldsymbol{\lambda}-2)}(L),
\end{equation}
defined with respect to the metric $g$.

We will restrict our attention to the two cases $\boldsymbol{\lambda}\in
(2-m,0)$ or $\boldsymbol{\lambda}\in (0,2)$. In either case $\mathcal{M}_L$ is
locally homeomorphic to the zero set of a smooth map 
$\Phi:\mathcal{I}\rightarrow\mathcal{O}$ defined (locally) between
finite-dimensional vector spaces. If furthermore $\boldsymbol{\mu}=2+\epsilon$
and all singularities are stable then $\mathcal{O}=\{0\}$ and $\mathcal{M}_L$ is
smooth of dimension $\mbox{dim}(\mathcal{I})$. Specifically: 
\begin{enumerate}
\item If $\boldsymbol{\lambda}\in (2-m,0)$ then
$\mbox{dim}(\mathcal{I})=b^1_c(L)-s$.
\item If $\boldsymbol{\lambda}\in (0,2)$ then
$\mbox{dim}(\mathcal{I})=b^1_{c,\bullet}(L)-s+\sum_{i=1}^l d_i$,
where $d_i$ is the number of harmonic functions on the AC end $S_i$ of the form
$r^\gamma\sigma(\theta)$ with $\gamma\in [0,\lambda_i]$.
\end{enumerate}
\end{theorem}
\begin{proof}
We split the proof into two parts, depending on the range of
$\boldsymbol{\lambda}$. To begin, assume $\boldsymbol{\lambda}\in (2-m,0)$. Using the notation of Decomposition \ref{decomp:closedforms_accs}, consider the (locally-defined) map
\begin{eqnarray*}
\tilde{F}:\tilde{\E}\times \widetilde{H}_{0,\infty}\times E_{0,\infty}\times
W^p_{k,(\boldsymbol{\mu},\boldsymbol{\lambda})}(L) &\rightarrow&
W^p_{k-2,(\boldsymbol{\mu}-2,\boldsymbol{\lambda}-2)}(L)\\
(\tilde{e},\beta,v,f) &\mapsto& F(\tilde{e},\beta+dv+df).
\end{eqnarray*}
Notice that $\tilde{F}$ is invariant under translations in $\R\subset E_{0,\infty}$. 
By regularity and Decomposition \ref{decomp:closedforms_accs}, $\mathcal{M}_L$ is
locally homeomorphic to $\tilde{F}^{-1}(0)/\R$. As in Proposition \ref{prop:accsnonlinear},
$d\tilde{F}[0](y,\beta,v,f)=d^*\beta+\Delta_g(\chi(y)+v+f)$. Now consider the
restricted map
\begin{equation}\label{eq:accsrestrictedlin}
d\tilde{F}[0]:T_e\tilde{\E}\oplus E_0\oplus
W^p_{k,(\boldsymbol{\mu},\boldsymbol{\lambda})}(L)\rightarrow
W^p_{k-2,(\boldsymbol{\mu}-2,\boldsymbol{\lambda}-2)}(L),
\end{equation}
where $E_0$ is the subspace of functions in $E_{0,\infty}$ which vanish on the
AC ends.  
We claim that this map is injective. To
prove this, assume $d\tilde{F}[0](y+v+f)=0$, \textit{i.e.} $\Delta_g(\chi(y)+v+f)=0$. Notice that $\chi(y)+v+f\in
W^p_{k,(-\boldsymbol{\epsilon},\boldsymbol{\lambda})}(L)$. Corollary \ref{cor:laplaceresults} shows
that $\Delta_g$ is an isomorphism on this space, so $\chi(y)+v+f=0$. In particular $d(\chi(y)+v+f)=0$ so
the infinitesimal Lagrangian deformation of $L$ defined by $(y,v,f)$ is trivial. This implies $y=0$ thus $\chi(y)=0$ and it is simple to conclude
that $f=0$ and $v=0$. 

Let $\mathcal{O}$ denote the cokernel of the map of Equation
\ref{eq:accsrestrictedlin}. More precisely, we define it to be a
finite-dimensional subspace of
$W^p_{k-2,(\boldsymbol{\mu}-2,\boldsymbol{\lambda}-2)}(L)$ such that 
\begin{equation}
\mathcal{O}\oplus d\tilde{F}[0]\left(T_e\tilde{\E}\oplus E_0\oplus
W^p_{k,(\boldsymbol{\mu},\boldsymbol{\lambda})}\right)=W^p_{k-2,(\boldsymbol{\mu
}-2,\boldsymbol{\lambda}-2)}(L).
\end{equation}
Consider the map
\begin{eqnarray*}
G:\mathcal{O}\times \tilde{\E}\times \widetilde{H}_{0,\infty}\times
E_{0,\infty}\times W^p_{k,(\boldsymbol{\mu},\boldsymbol{\lambda})}(L)
&\rightarrow& W^p_{k-2,(\boldsymbol{\mu}-2,\boldsymbol{\lambda}-2)}(L)\\
(\gamma,\tilde{e},\beta,v,f) &\mapsto& \gamma+\tilde{F}(\tilde{e},\beta,v,f).
\end{eqnarray*}
Again, $G$ is invariant under translations in $\R$.
By construction the restriction of $dG[0]$ to the space $\mathcal{O}\oplus
T_e\tilde{\E}\oplus E_0\oplus W^p_{k,(\boldsymbol{\mu},\boldsymbol{\lambda})}$
is an isomorphism. We now have
the following information about the map $G$. First,
let $E'$ denote a complement of $E_0\oplus\R$ in
$E_{0,\infty}$, \textit{i.e.} $E_{0,\infty}=E_0\oplus\R\oplus E'$. Then $\mbox{Ker}(dG[0])=V\oplus\R$, where $V$ is some vector space projecting
isomorphically onto $\widetilde{H}_{0,\infty}\oplus E'$. Second, by the Implicit Function
Theorem, the set $G^{-1}(0)$ is smooth and can be locally written as the graph
of a smooth map $\Phi$ defined on the kernel of $dG[0]$, thus on $\widetilde{H}_{0,\infty}\oplus(\R\oplus E')$.
As in Proposition \ref{prop:fredholmreducestofinitedim} we can conclude that the
projection onto $\widetilde{H}_{0,\infty}\oplus(\R\oplus E')$ restricts to a homeomorphism
$\tilde{F}^{-1}(0)\simeq(\pi_{\mathcal{O}}\circ\Phi)^{-1}(0)$. It is simple to
check that $\Phi$ is invariant under translations in $\R$.
Restricting $\Phi$ to $\mathcal{I}:=\widetilde{H}_{0,\infty}\oplus E'$ proves the first
claim regarding $\mathcal{M}_L$ for this range of
$\boldsymbol{\lambda}$. Notice that
$\mbox{dim}(\widetilde{H}_{0,\infty})=b^1_c(L)-(s+l)+1$ and $\mbox{dim}(E')=l-1$
so $\mbox{dim}(\mathcal{I})=b^1_c(L)-s$.

Now let us further assume that $\boldsymbol{\mu}=2+\epsilon$ and that all
singularities are stable. Here $\epsilon$ is to
be understood as in Remark \ref{rem:muregularity}. By Corollary
\ref{cor:laplaceresults} and the definition of stability,
\begin{equation}
\mbox{dim(Coker$(\Delta_{\boldsymbol{\mu},\boldsymbol{\lambda}})$)}=d, \ \
\mbox{where } d:=\sum_{i=1}^s \left(1+2m+m^2-1-\mbox{dim}(G_i)\right).
\end{equation}
Again, $d$ is also the dimension of the space $T_e\tilde{\E}\oplus E_0$. Our
previous injectivity calculation thus implies that the map $d\tilde{F}[0]$ of
Equation \ref{eq:accsrestrictedlin} is an isomorphism. In particular,
$\mathcal{O}=\{0\}$. We can now apply the Implicit Function Theorem directly to
$\tilde{F}$ to obtain that $\tilde{F}^{-1}(0)$ is smooth. Quotienting by $\R$
shows that $\mathcal{M}_L$ is smooth.

We now start over again, under the assumption $\boldsymbol{\lambda}\in (0,2)$.
In this case we use the map
\begin{eqnarray*}
\tilde{F}:\tilde{\E}\times \widetilde{H}_{0,\bullet}\times E_0\times
W^p_{k,(\boldsymbol{\mu},\boldsymbol{\lambda})}(L) &\rightarrow&
W^p_{k-2,(\boldsymbol{\mu}-2,\boldsymbol{\lambda}-2)}(L)\\
(\tilde{e},\beta,v,f) &\mapsto& F(\tilde{e},\beta+dv+df)
\end{eqnarray*}
and the restricted map
\begin{equation}\label{eq:accsrestrictedlinbis}
d\tilde{F}[0]:T_e\tilde{\E}\oplus E_0\oplus
W^p_{k,(\boldsymbol{\mu},\boldsymbol{\lambda})}(L)\rightarrow
W^p_{k-2,(\boldsymbol{\mu}-2,\boldsymbol{\lambda}-2)}(L).
\end{equation}
Recall the construction of $E_0$ in Decomposition \ref{decomp:closedforms_accs}:
it is clear that we may assume that $\chi(T_e\tilde{\E})$ and $E_0$ are linearly
independent in $W^p_{k,(-\boldsymbol{\epsilon},-\boldsymbol{\epsilon})}(L)$.
Corollary \ref{cor:laplaceresults} proves that $\Delta_g$ is injective on this
space. Define a decomposition
\begin{equation}
T_e\tilde{\E}\oplus E_0=Z'\oplus Z''
\end{equation}
by imposing $\Delta_g(Z')=\Delta_g(T_e\tilde{\E}\oplus E_0)\cap
\mbox{Im}(\Delta_{\boldsymbol{\mu},\boldsymbol{\lambda}})$ and choosing any
complement $Z''$. Then one can check that the kernel of the map of Equation
\ref{eq:accsrestrictedlinbis} is isomorphic to
$Z'\oplus\mbox{Ker}(\Delta_{\boldsymbol{\mu},\boldsymbol{\lambda}})$.

Choose $\mathcal{O}$ in
$W^p_{k-2,(\boldsymbol{\mu}-2,\boldsymbol{\lambda}-2)}(L)$ such that 
\begin{equation}
\mathcal{O}\oplus d\tilde{F}[0]\left(T_e\tilde{\E}\oplus E_0\oplus
W^p_{k,(\boldsymbol{\mu},\boldsymbol{\lambda})}\right)=W^p_{k-2,(\boldsymbol{\mu
}-2,\boldsymbol{\lambda}-2)}(L).
\end{equation}
Consider the map
\begin{eqnarray*}
G:\mathcal{O}\times \tilde{\E}\times \widetilde{H}_{0,\bullet}\times E_0\times
W^p_{k,(\boldsymbol{\mu},\boldsymbol{\lambda})}(L) &\rightarrow&
W^p_{k-2,(\boldsymbol{\mu}-2,\boldsymbol{\lambda}-2)}(L)\\
(\gamma,\tilde{e},\beta,v,f) &\mapsto& \gamma+\tilde{F}(\tilde{e},\beta,v,f).
\end{eqnarray*}
The restriction of $dG[0]$ to the space $\mathcal{O}\oplus T_e\tilde{\E}\oplus
E_0\oplus W^p_{k,(\boldsymbol{\mu},\boldsymbol{\lambda})}$ is surjective. As
before, this implies that $G^{-1}(0)$ can be parametrised via a smooth map
$\Phi$ defined (locally) on the space $\widetilde{H}_{0,\bullet}\oplus Z'\oplus
\mbox{Ker}(\Delta_{\boldsymbol{\mu},\boldsymbol{\lambda}})$. As usual, these
maps are invariant under translations in $\R\subset
Z'\oplus\mbox{Ker}(\Delta_{\boldsymbol{\mu},\boldsymbol{\lambda}})$. Setting
$\mathcal{I}:=(\widetilde{H}_{0,\bullet}\oplus Z'\oplus
\mbox{Ker}(\Delta_{\boldsymbol{\mu},\boldsymbol{\lambda}}))/\R$ and considering
the natural map on this quotient then proves the first claim regarding
$\mathcal{M}_L$ for this range of $\boldsymbol{\lambda}$.

Now assume that $\boldsymbol{\mu}=2+\epsilon$ and that all singularities are
stable. Choose $\boldsymbol{\lambda}'\in (2-m,0)$. We can restrict the map of
Equation \ref{eq:accsrestrictedlinbis} to the map
\begin{equation}\label{eq:accsrestrictedlinter}
d\tilde{F}[0]:T_e\tilde{\E}\oplus E_0\oplus
W^p_{k,(\boldsymbol{\mu},\boldsymbol{\lambda}')}(L)\rightarrow
W^p_{k-2,(\boldsymbol{\mu}-2,\boldsymbol{\lambda}'-2)}(L).
\end{equation}
Exactly as for Equation \ref{eq:accsrestrictedlin}, it is simple to prove that
Equation \ref{eq:accsrestrictedlinter} defines an isomorphism and that
$\mbox{dim}(T_e\tilde{\E}\oplus
E_0)=\mbox{dim(Coker($\Delta_{\boldsymbol{\mu},\boldsymbol{\lambda}'}$))}$,
where 
\begin{equation*}
\Delta_{\boldsymbol{\mu},\boldsymbol{\lambda}'}:=\Delta_g:
W^p_{k,(\boldsymbol{\mu},\boldsymbol{\lambda}')}(L)\rightarrow
W^p_{k-2,(\boldsymbol{\mu}-2,\boldsymbol{\lambda}'-2)}(L).
\end{equation*}
One can check that the dimension of
$\mbox{Coker}(\Delta_{\boldsymbol{\mu},\boldsymbol{\lambda}})$ decreases as
$\boldsymbol{\lambda}$ increases. We can actually assume, cf.
\cite{pacini:weighted}, that
$\mbox{Coker}(\Delta_{\boldsymbol{\mu},\boldsymbol{\lambda}})\subseteq\mbox{
Coker}(\Delta_{\boldsymbol{\mu},\boldsymbol{\lambda}'})$. This proves that the
map of Equation \ref{eq:accsrestrictedlinbis} is surjective, \textit{i.e.}
$\mathcal{O}=\{0\}$, so $\tilde{F}^{-1}(0)$ and $\mathcal{M}_L$ are smooth. To
compute the dimension of this moduli space notice that $Z''\simeq
\mbox{Coker}(\Delta_{\boldsymbol{\mu},\boldsymbol{\lambda}})$ so
\begin{eqnarray}
\mbox{dim(Ker$(d\tilde{F}[0])$)} &=&
\mbox{dim(Ker$(\Delta_{\boldsymbol{\mu},\boldsymbol{\lambda}})$)}
+\mbox{dim}(Z')\nonumber\\
&=&
\mbox{dim(Ker$(\Delta_{\boldsymbol{\mu},\boldsymbol{\lambda}})$)}+\mbox{
dim(Coker($\Delta_{\boldsymbol{\mu},\boldsymbol{\lambda}'}$))}-\mbox{
dim(Coker($\Delta_{\boldsymbol{\mu},\boldsymbol{\lambda}}$))}\nonumber\\
&=&
i(\Delta_{\boldsymbol{\mu},\boldsymbol{\lambda}})-i(\Delta_{\boldsymbol{\mu},
\boldsymbol{\lambda}'}),
\end{eqnarray}
where $i$ denotes the index of the Fredholm map. This implies that the kernel of
the full map $d\tilde{F}[0]$ has dimension
$\mbox{dim}(\widetilde{H}_{0,\bullet})+i(\Delta_{\boldsymbol{\mu},\boldsymbol{
\lambda}})-i(\Delta_{\boldsymbol{\mu},\boldsymbol{\lambda}'})$. The conclusion
follows from Equation \ref{eq:cohomsequence_accsbis} and the change of index
formula described in \cite{pacini:weighted}.
\end{proof}
We call $\mathcal{O}$ the \textit{obstruction space} of the SL deformation
problem.

\begin{example}
Let $\mathcal{C}$ be a SL cone in $\C^m$. Assume $\mathcal{C}$ is stable and
that its link $\Sigma$ is connected so that $s=1$. Using Poincar\'{e} Duality
and the fact that $\mathcal{C}\simeq\Sigma\times (0,\infty)$ we see that 
\begin{equation}
b^1_c(\mathcal{C})=b^{m-1}(\mathcal{C})=b^{m-1}(\Sigma)=1.
\end{equation}
Theorem \ref{thm:accssl} then shows that, for $\lambda\in (2-m,0)$,
$\mathcal{M}_{\mathcal{C}}$ has dimension 0, \textit{i.e.} $\mathcal{C}$ is
rigid within this class of deformations.

Notice also that restriction defines isomorphisms $H^i(\mathcal{C};\R)\simeq
H^i(\Sigma;\R)$ so the long exact sequence \ref{eq:cohomsequence_accsbis}, using
$\Sigma_0=\Sigma$, leads to $H^i_{c,\bullet}(\mathcal{C};\R)=0$. 
Theorem \ref{thm:accssl} then shows that $\mathcal{M}_{\mathcal{C}}$ has
dimension 0 if $\lambda\in (0,1)$ and has dimension $2m$ if $\lambda\in (1,2)$.
In the latter case the SL deformations are simply the translations of
$\mathcal{C}$ in $\C^m$.
\end{example}


\subsection{Comparison to other results in the literature}\label{ss:comparison}
It is interesting to compare Theorem \ref{thm:accssl} to other moduli space results available in the literature.
The first such result, for compact SLs in a CY manifold $M$, was proved by McLean \cite{mclean}.

\begin{theorem} \label{thm:mclean}
Let $L$ be a smooth compact SL submanifold of a CY manifold $M$. Let
$\mathcal{M}_L$ denote the moduli space of SL deformations of $L$. Then
$\mathcal{M}_L$ is a smooth manifold of dimension $b^1(L)$. 
\end{theorem}
A special feature of this compact setting is the fact that $dF[0]=d^*$ is not surjective. In theory this should interfere with the Implicit Function Theorem argument. However this is actually not a problem because Proposition \ref{prop:cptnonlinear} (2) allows us to also restrict the range of $F$; the linearization of the restricted map is surjective. Theorem \ref{thm:mclean} can thus be proved similarly to Theorem \ref{thm:accssl}. We define the map $\tilde{F}$ using $K:=H$, the space determined in Decomposition \ref{decomp:closedforms}. Restricted to functions, \textit{i.e.} to the Hamiltonian deformations of L, $d\tilde{F}[0]=\Delta_g$ is an isomorphism (after restricting the range of $\tilde{F}$ as above), cf. Theorem \ref{thm:laplaceresults}, so $\mathcal{M}_L\simeq\tilde{F}^{-1}(0)$ can be written as a smooth graph over $K$, proving the result and the dimension count $b^1(L)$.

The corresponding result for AC SLs was proved independently by the
author \cite{pacini:defs} and by Marshall \cite{marshall}. 

\begin{theorem} \label{thm:pacini}
Let $L$ be an AC SL submanifold of $\C^m$ with rate
$\boldsymbol{\lambda}$. Let $\mathcal{M}_L$ denote the moduli space of SL
deformations of $L$ with rate $\boldsymbol{\lambda}$. Consider the operator
\begin{equation}\label{eq:aclaplacian}
\Delta_g:W^p_{k,\boldsymbol{\lambda}}(L)\rightarrow
W^p_{k-2,\boldsymbol{\lambda}-2}(L).
\end{equation}
\begin{enumerate}
\item If $\boldsymbol{\lambda}\in (0,2)$ is a non-exceptional weight for
$\Delta_g$ then $\mathcal{M}_L$ is a smooth manifold of dimension
$b^1(L)+\mbox{dim(Ker$(\Delta_g)$)}-1$. 
\item If $\boldsymbol{\lambda}\in (2-m,0)$ then $\mathcal{M}_L$ is a smooth
manifold of dimension $b^1_c(L)$.
\end{enumerate}
\end{theorem}
This result can be obtained as a special case of Theorem \ref{thm:accssl} simply assuming that the set of CS ends is empty, \textit{i.e.} $s=0$. When $\boldsymbol{\lambda}<0$ we may define $\tilde{F}$ using $K:=\tilde{H}_{\infty}\times d(E_{\infty})$, cf. Decomposition \ref{decomp:closedforms_decay}. Restricted to the complement of $K$, Theorem \ref{thm:laplaceresults} shows that $d\tilde{F}[0]$ is an isomorphism so $\mathcal{M}_L$ is parametrized by $K$. When $\boldsymbol{\lambda}\in (0,2)$ we set $K:=H$, cf. Decomposition \ref{decomp:closedforms_growth}. Restricted to the complement of $K$, Theorem \ref{thm:laplaceresults} shows that $d\tilde{F}[0]$ is surjective but it has kernel which contributes to the parameters defining $\mathcal{M}_L$. In both cases $K$ corresponds exactly to the ``topological'' contributions to the dimension count, as emphasized in Corollary \ref{cor:topsummary}. It is interesting to notice however that in the case $\boldsymbol{\lambda}<0$ the space $K$ contains some Hamiltonian contributions,  corresponding to $d(E_{\infty})$.

Finally, Joyce \cite{joyce:II} proved the following result on CS SLs in general CYs.
\begin{theorem} \label{thm:joyce}
Let $L$ be a CS SL submanifold of $M$ with $s$ singularities and rate
$\boldsymbol{\mu}$. Let $\mathcal{M}_L$ denote the moduli space of SL
deformations of $L$ with moving singularities and rate $\boldsymbol{\mu}$.
Assume $\boldsymbol{\mu}$ is non-exceptional for the map
\begin{equation}\label{eq:cslaplacian}
\Delta_g:W^p_{k,\boldsymbol{\mu}}(L)\rightarrow \{u\in
W^p_{k-2,\boldsymbol{\mu}-2}(L): \int_L u\,vol_g=0\}.
\end{equation}
Then $\mathcal{M}_L$ is locally homeomorphic to the zero set of a smooth map
$\Phi: \mathcal{I}\rightarrow\mathcal{O}$
defined (locally) between finite-dimensional vector spaces. If
$\boldsymbol{\mu}=2+\epsilon$ and all singularities are stable then
$\mathcal{O}=\{0\}$ and $\mathcal{M}_L$ is smooth of dimension
$\mbox{dim}(\mathcal{I})=b^1_c(L)-s+1$.
\end{theorem}
In this case we can set $K:=\tilde{\E}\times \tilde{H}_0\times d(E_0)$ (cf. Decomposition \ref{decomp:closedforms_decay}). Then the stability condition implies that, after restricting the range of $\tilde{F}$ as in the smooth compact case, $d\tilde{F}[0]$ is an isomorphism on $T_e\tilde{\E}\oplus d(E_0)\oplus
W^p_{k,\boldsymbol{\mu}}(L)$ so $\mathcal{M}_L$ is parametrized by $\tilde{H}_0$, whose dimension is calculated in Corollary \ref{cor:topsummary}.

\subsubsection*{Concluding remarks} When $\boldsymbol{\lambda}<0$ and the stability condition is
verified, the dimension of the SL moduli spaces appearing in Theorems
\ref{thm:accssl}, \ref{thm:mclean}, \ref{thm:pacini} and \ref{thm:joyce} is
purely topological. The cases analyzed in the theorems correspond exactly to the
cases analyzed in Corollary \ref{cor:topsummary}, in the sense that the moduli
spaces should be thought of as being modelled on the cohomology spaces which
appear in Corollary \ref{cor:topsummary}. 

It is interesting to notice how decay conditions on AC and CS ends are
incorporated differently into these cohomology spaces: decay conditions on AC
ends correspond to using compactly-supported forms while decay conditions on CS
ends correspond to the condition that a certain restriction map vanishes, cf. also Remark \ref{rem:topsummary}.

Allowing $\boldsymbol{\lambda}>0$ changes the topological data, again in
agreement with Corollary \ref{cor:topsummary}. It also introduces new SL
deformations which depend on analytic data.

Finally we point out that the role of the space $\tilde{\E}$ (thus of the stability condition) is always to contribute to making the linearized operator surjective. This means that $\tilde{\E}$ never contributes parameters to the moduli space. In other words, the position of the singularities and the direction of the CS cones of the deformed submanifolds are forced by the analysis, and cannot be assigned arbitrarily. Translations of the AC ends correspond instead to harmonic functions of linear growth, so they appear among  the parameters of the moduli space as soon as $\boldsymbol{\lambda}>1$.


\ 

{\bf Acknowledgments.} I would like to thank D. Joyce for many useful
explanations on his work and for his help and advice on various parts of this
paper. I would also like to thank S. Karigiannis and Y. Song for useful conversations. The main part of this 
work was carried out while I was a Marie Curie EIF Fellow at the University of
Oxford.

\bibliographystyle{amsplain}
\bibliography{accssldefs}
\end{document}